\documentclass{article}
\usepackage{hyperref}
\usepackage{epsfig}
\usepackage{amsmath}
\usepackage{amsfonts}
\usepackage{amsthm}
\usepackage{algorithm}
\usepackage{algorithmic}
\usepackage{graphicx}
\usepackage{subfig}
\usepackage{array}
\usepackage{rotating}
\usepackage[T1]{fontenc}

\usepackage{pgfplots,tikz}
\usetikzlibrary{plotmarks}

\newcommand {\C}        {{\mathbb{C}}}
\newcommand {\R}        {{\mathbb{R}}}
\newcommand {\N}        {{\mathbb{N}}}

\newcommand {\cH}  {\mathcal{H}}
\newcommand {\cL}  {\mathcal{L}}
\newcommand {\ri}  {\mathrm{i}}

\newcommand{\abs}[1]{\left|#1\right|}

\DeclareMathOperator{\Col}{Col}

\DeclareMathOperator{\Real}{Re}

\newtheorem{theorem}{Theorem}[section]
\newtheorem{lemma}[theorem]{Lemma}

\newtheorem{example}[theorem]{Example}
\newtheorem{remark}[theorem]{Remark}

\title{Large-Scale Computation of $\cL_\infty$-Norms by a Greedy Subspace Method}
\author{
Nicat Aliyev\footnote{Ko\c{c} University, Department of Mathematics, Rumeli Feneri Yolu 34450, Sar{\i}yer, Istanbul, Turkey, E-Mail: \texttt{naliyev@ku.edu.tr}.} \and
Peter Benner\footnote{Max Planck Institute for Dynamics of Complex Technical Systems, Sandtorstra{\ss}e 1, 39106 Magdeburg, Germany, E-Mail: \texttt{benner@mpi-magdeburg.mpg.de}.} \and
Emre Mengi\footnote{Ko\c{c} University, Department of Mathematics, Rumeli Feneri Yolu 34450, Sar{\i}yer, Istanbul, Turkey, E-Mail: \texttt{emengi@ku.edu.tr}.
The work of the author was supported in part by the BAGEP program of Turkish Academy of Science.} \and
Paul Schwerdtner\footnote{Technische Universit\"at Berlin, Institut f\"ur Mathematik, Stra{\ss}e des 17. Juni 136, 10623 Berlin, Germany, E-Mail: \texttt{p.schwerdtner@campus.tu-berlin.de}. This work is supported by the DFG priority program 1897: ``Calm, Smooth and Smart -- Novel Approaches for Influencing Vibrations by Means of Deliberately Introduced Dissipation''.} \and
Matthias Voigt\footnote{Technische Universit\"at Berlin, Institut f\"ur Mathematik, Stra{\ss}e des 17. Juni 136, 10623 Berlin, Germany, E-Mail: \texttt{mvoigt@math.tu-berlin.de}. This work is supported by the Einstein Foundation Berlin within the framework of the Einstein Center for Mathematics (ECMath).}
}

\begin{document}
\maketitle
 

\begin{abstract}
\noindent
We are concerned with the computation of the $\cL_\infty$-norm for an $\cL_\infty$-function of the form
$H(s) = C(s) D(s)^{-1} B(s)$, where the middle factor is the inverse of a meromorphic matrix-valued function, 
and $C(s),\, B(s)$ are meromorphic functions mapping to short-and-fat and tall-and-skinny matrices, respectively. 
For instance, transfer functions of descriptor systems and delay systems fall into this family. We focus on 
the case where the middle factor is {large-scale}. We propose a subspace projection method to obtain 
approximations of the function $H$ where the middle factor is of much smaller dimension. The $\cL_\infty$-norms 
are computed for the resulting reduced functions, then the subspaces are refined by means of the optimal 
points on the imaginary axis where the $\cL_\infty$-norm of the reduced function is attained. 
The subspace method is designed so that certain Hermite interpolation properties hold between the largest 
singular values of the original and reduced functions. This leads to a locally superlinearly convergent algorithm with 
respect to the subspace dimension, which we prove and illustrate on various numerical examples.

\vskip 1ex



\vskip 1ex

\end{abstract}

\section{Introduction}
We consider the computation of the ${\mathcal L}_\infty$-norm of a matrix-valued function of the form
\begin{equation}\label{eq:transfer_func}
	H : \Omega \rightarrow {\mathbb C}^{p\times m},	\quad	H(s) :=  C(s) D(s)^{-1} B(s),
\end{equation}
and specifically address the case when the middle square factor {$D$ is of large dimension}. In what follows{,}
a subspace method is derived to reduce the size of {$D$} making efficient computation
of the ${\mathcal L}_\infty$-norm of {$H$} possible. The domain $\Omega$ is 
{an open} subset of the complex plane and assumed to enclose the imaginary axis $\ri\R$. 


Furthermore, it is assumed that the functions 
$B: \Omega \rightarrow {\mathbb C}^{n\times m}$, 
$C: \Omega \rightarrow {\mathbb C}^{p\times n}$, and 
$D : \Omega \rightarrow {\mathbb C}^{n\times n}$ {are 
defined by}
\begin{align}\label{eq:matsum}
 \begin{split}
	 B(s)  &:=  f_1(s) B_1 + \dots + f_{\kappa_B}(s) B_{\kappa_B}, \\
	 C(s)  &:=  g_1(s) C_1 + \dots + g_{\kappa_C}(s) C_{\kappa_C}, \\
	 D(s)  &:=  h_1(s) D_1 + \dots + h_{\kappa_D}(s) D_{\kappa_D},
 \end{split}
\end{align}
for given matrices $B_1,\, \dots,\, B_{\kappa_B} \in {\mathbb C}^{n\times m}$, $C_1,\, \dots,\, C_{\kappa_C} \in {\mathbb C}^{p\times n}$, $D_1,\, \dots,\, D_{\kappa_D} \in {\mathbb C}^{n\times n}$ and given functions $f_1,\, \dots, f_{\kappa_B},\,g_1,\, \dots,\, g_{\kappa_C},\,h_1,\,\dots,\, h_{\kappa_D} : \Omega \rightarrow {\mathbb C}$ that are assumed to be meromorphic in $\Omega$. 

For example, if $sE-A$ is a regular pencil, then the transfer function
\begin{equation*}
 H(s) = C(sE-A)^{-1}B
\end{equation*}
of the descriptor system
\begin{equation}\label{eq:disc_sys}
    E x'(t) = A x(t) + B u(t), \quad y(t) = C x(t),
\end{equation}
and more generally, the transfer function
\begin{equation*}
 H(s) = C\left( sE - A_0 - \sum_{j=1}^m \mathrm{e}^{-s\tau_j} A_j \right)^{-1} B
\end{equation*}
of the 
delay differential-algebraic system
\begin{equation}\label{eq:delay_sys}
	Ex'(t) = A_0 x(t) + \sum_{j=1}^m A_j x(t - \tau_j) + B u(t), \quad y(t) = C x(t)
\end{equation}
are encompassed by framework {\eqref{eq:transfer_func}--\eqref{eq:matsum}}.
Some other examples are transfer functions of 
higher order systems and systems containing input and output delays,
as well as transfer functions of the form 
\begin{equation*}
 H(s) =  sB^\ast (s^2I_n - \sqrt{s} D_2 + D_3)^{-1}B
\end{equation*}
resulting from the spatial discretization of
electromagnetic field equations, {i.\,e.,} the Maxwell equations,
describing the electro-dynamical behavior of microwave devices with
surface losses (see \cite{morFenB10} and references therein). 

We are concerned with the computation of the $\cL_\infty$-norm of $H$, particularly for the case where $n$ is very large and further $n \gg m,\,p$. 
We define the spaces
\begin{align*}
 \cL_\infty^{p \times m} &:= {\left\{ H|_{\mathrm{i}\R} \; \bigg| \; H : \Omega \rightarrow \C^{p \times m} \text{ is analytic for an open domain } \Omega \subseteq \C \text{ with } \ri\R \subset \Omega \right.} \\ & \hphantom{:= \left\{ H|_{\mathrm{i}\R} \; \bigg| \;\right.} \qquad  {\left. \text{ and } \sup_{\omega \in \R}\left\| H(\ri\omega) \right\|_2 < \infty \right\}}, \\
 \cH_\infty^{p \times m} &:= \left\{ H : \C^+ \rightarrow {\mathbb C}^{p\times m} \; \bigg| \; H \text{ is analytic and } \sup_{s \in \C^+} \left\| H(s) \right\|_2 < \infty \right\},
\end{align*}
{where $\C^+:=\{ s \in \C \;|\; \Real(s) > 0 \}$ denotes the open right complex half-plane.}
In this paper, {the function $H$ (more precisely, its restriction to the imaginary axis)} is assumed to be in $\cL_\infty^{p \times m}$. {For ease of notation, we write $H \in \cL_\infty^{p \times m}$ instead of $H|_{\ri\R} \in \cL_\infty^{p \times m}$.}
For such, the $\cL_\infty$-norm is defined by
\begin{equation*}
 \left\| H \right\|_{\cL_\infty} := \sup_{\omega \in \R} \left\| H(\ri\omega) \right\|_2 = \sup_{\omega \in \R} \sigma(H(\ri \omega)), 
\end{equation*}
where $\sigma(\cdot)$ denotes the largest singular value of its matrix argument.
Throughout the text we refer each function in $\cL_\infty^{p \times m}$ as an \emph{$\cL_\infty$-function}.

In most applications one is often rather interested in functions which are in $\cH_\infty^{p \times m}$. For such, using the maximum principle for analytic functions, one can show that the $\cH_\infty$-norm is equivalent to the $\cL_\infty$-norm, i.\,e.,
\begin{equation*}
 \left\| H \right\|_{\cH_\infty} := \sup_{s \in \C^+} \left\| H(s) \right\|_2 = \sup_{s \in \partial\C^+} \left\| H(s) \right\|_2 = \sup_{\omega \in \R} \sigma(H(\ri \omega)). 
\end{equation*}

\subsection{Motivation}
The $\cH_\infty$-norm plays an indispensable role in the assessment of robust stability as well as in robust control. For instance, assume that we are given an exponentially stable delay differential-algebraic equation
\begin{equation*}
 E x'(t) = A_0 x(t) + A_1 x(t - \tau),
\end{equation*}
and consider the perturbed delay differential-algebraic equation \cite{DuLMT13}
\begin{equation} \label{eq:pertDDAE}
 (E + B_1 \Delta_1 C) x'(t) = (A_0 + B_2 \Delta_2 C) x(t) + (A_1 + B_3 \Delta_3 C) x(t - \tau),
\end{equation}
where $\Delta_i \in \C^{m_i \times p}$, $i = 1,\,2,\,3$ are the perturbations and $B_i \in \C^{n \times m_i}$ and $C \in \C^{p \times n}$ are matrices that define the perturbation structure. Define the function 
\begin{equation*}
 H(s) := C (sE - A_0 - \mathrm{e}^{-s\tau} A_1)^{-1} \begin{bmatrix} -sB_1 & B_2 & \mathrm{e}^{-s \tau}B_3 \end{bmatrix}.
\end{equation*}
In \cite{DuLMT13} it is shown that under certain conditions on the
matrices $E$, $A_0$, $A_1$ (ensuring a ``strangeness-free'' system) and
some further restrictions on the perturbation structure matrices
$B_1$, $B_2$, $B_3$, the $\cH_\infty$-norm is the reciprocal of the
structured stability radius, similar to the standard state-space
case \cite{HinP86a,HinP86b}. In other words, with $\Delta
:= \begin{bmatrix} \Delta_1^* & \Delta_2^* &
  \Delta_3^* \end{bmatrix}^*$ we have 
\begin{equation*}
 \inf\left\{ \left\| \Delta \right\|_2 \; \big| \; \text{system \eqref{eq:pertDDAE} is not exponentially stable} \right\} = \left\| H \right\|_{\cH_\infty}^{-1}.
\end{equation*}

This connection also motivates the importance of the $\cH_\infty$-norm in robust control, and
the minimization of the $\cH_\infty$-norm over system parameters. Consider, for example, the system (see, e.\,g., \cite{ZhoDG96}),
\begin{align*}
  E x'(t) &= A_0x(t) + B_1u(t) + B_2w(t), \\
     y(t) &= C_1x(t-\tau), \\
     z(t) &= C_2x(t),
\end{align*}
where $u$ is the control input, $y$ is the (delayed) measured output, $w$ is an input representing noise or unmodeled dynamics, and $z$ is the performance output, respectively. By imposing the feedback law $u(t) = Fy(t)$, we obtain the closed-loop system
\begin{align*}
  E x'(t) &= A_0x(t) + B_1FC_1x(t-\tau) + B_2w(t), \quad
     z(t) = C_2x(t).
\end{align*}
With $A_{1,F} := B_1FC_1$ its transfer function from $w$ to $z$ is given by
\begin{equation*}
 H_F(s) = C_2(sE-A_0-\mathrm{e}^{-s\tau}A_{1,F})^{-1}B_2.
\end{equation*}
The goal of robust control is to determine a stabilizing feedback $F$ such that the closed-loop $\cH_{\infty}$-norm, i.\,e., $\left\| H_F \right\|_{\cH_\infty}$ is minimized in order to achieve a maximum robustness of stability of the performance output $z$ with respect to disturbances and noise that enter the system via the input $w$. For standard state-space systems this $\cH_\infty$ optimization problem is addressed by the MATLAB package \textsc{HIFOO} \cite{BurHLO06}. In the past ten years this software has found manifold applications in industry, some of which are outlined in \cite{MitO15}. Since \textsc{HIFOO} performs a couple of $\cH_\infty$-norm evaluations, an efficient $\cH_\infty$-norm computation will be beneficial for the performance of the optimization procedure. 

\subsection{Literature}
Studies concerning the computation of the $\cL_\infty$- or $\cH_\infty$-norm have been
conducted since the late 1980s. Byers' work \cite{Byers1988} focuses
on the computation of the distance to instability for a matrix, which
can be viewed as a special $\cL_\infty$-norm computation problem for the transfer function of a
standard state-space system \eqref{eq:disc_sys} with $B$, $C$, and $E$
being 
identities. This idea has been independently adapted for the computation
of the $\cL_\infty$-norm of transfer functions of standard state-space systems by Boyd,
Balakrishnan \cite{Boyd1990}, as well as Bruinsma and Steinbuch
\cite{Bruinsma1990}. An extension of these methods to transfer functions of descriptor
systems is discussed in \cite{BenSV12}. These are level-set based
optimization approaches, and require the repeated solution of
Hamiltonian eigenvalue problems of size twice the order of the
system. Consequently, they are not suitable for systems beyond medium
scale.  

For larger problems, several approaches have been proposed in recent years. For instance, the characterization of the $\cL_\infty$-norm via a Hamiltonian eigenvalue problem has been used to formulate an associated root-finding problem which can be solved using Newton's method \cite{Freitag2014}. This approach requires solutions of linear systems of size equal to the order of the system.
Some other approaches \cite{Guglielmi2013, Voigt2015, MitO15b} are restricted to the case of the $\cH_\infty$-norm only. They are based on the relation of the $\cH_\infty$-norm to the structured stability radius and structured 
$\varepsilon$-pseudospectra \cite{Trefethen2005, Karow2003}; these approaches compute the rightmost point of the structured $\varepsilon$-pseudospectrum repeatedly for various values of $\varepsilon$. However, all of these methods for larger problems converge only locally and there is no guarantee that the global maximum of $\sigma(H(\ri \cdot))$ is found.  

The delay-system setting is addressed by a few works \cite{Gumussoy2009, Gumussoy2015} only. These are extensions of
the level-set based approach of Byers, but involve infinite dimensional operators. None of these works benefits from a subspace projection idea and their use is typically limited to systems of the order of thousand at most.

\subsection{Contributions and Outline}
Our approach is based on a reduction of the middle factor $D(s)$ in \eqref{eq:transfer_func} to a much smaller dimension using two-sided projections.
The $\cL_\infty$-norm is computed for the resulting reduced matrix-valued function, then the subspaces are expanded using the singular vectors of $H(\ri\omega_r)$, where $\ri\omega_r$ is the point on the imaginary axis (including infinity) at which the reduced function attains its $\cL_\infty$-norm.  
Our expansion strategy leads to superlinear convergence with respect to the subspace dimension which we observe in practice and prove in theory. 
This work is inspired by a recent work \cite{Kangal2015} on a subspace method in the context of eigenvalue optimization. 
However, unlike \cite{Kangal2015}, the matrix-valued function $H(\ri \cdot)$ (whose largest singular value is to be maximized)
is $p\times m$ where $p,\,m$ are typically small, the large-scale nature of the problem in this paper is due to the size of $D(s)$. Dealing with the large dimensionality of $D(s)$
requires a different approach compared to the one proposed for eigenvalue optimization in \cite{Kangal2015}.

We expose our work in the following order. In the next section, we formally introduce the reduced matrix-valued functions and present a result (Theorem~\ref{thm:transfer_interpolate}) that points out how Hermite interpolation of the original $\cL_\infty$-function can be achieved by a reduced 
matrix-valued function. This interpolation result gives rise to the formal definition of the subspace method as in Algorithm~\ref{alg1}. The method is devised in order to lift the Hermite interpolation properties to the largest singular value functions associated with the original $\cL_\infty$-function and the reduced matrix-valued function. The local superlinear convergence of the subspace method can be attributed to these interpolation properties. This convergence is proven rigorously in Section \ref{sec:convergence}. Important implementation details of the proposed method and the results of our numerical experiments are discussed in Section \ref{sec:experiment}.

\section{Our Approach}\label{sec:approach}
Two-sided subspace projections are widely used in model order reduction \cite{DeVillemagne1987, Yousuff1985, morBauBF14, Gugercin2013}.
In the context of a descriptor system of the form \eqref{eq:disc_sys}, this amounts to restricting the state-space to a subspace ${\mathcal V}$ of dimension much smaller than the original state-space, and imposing a Petrov-Galerkin condition
with respect to another subspace ${\mathcal W}$. 
Formally, introducing matrices $V,\, W$ whose columns span ${\mathcal V},\, {\mathcal W}$, respectively, 
the reduced state at time $t$ is given by $V \tilde{x}(t)$, and the reduced system is defined by
\[
	W^\ast 
	\left(
		EV \tilde{x}'(t)	- AV \tilde{x}(t)  - B u(t)	
	\right)
			=
	0 \quad	{\rm and}	\quad
	y(t) = C V \tilde{x}(t).
\]
The transfer functions associated with the original descriptor system and the reduced one above are
\[
	H(s) = C(sE - A)^{-1}B \quad	{\rm and} \quad
	\widetilde{H}(s) = CV (sW^\ast E V - W^\ast AV )^{-1} W^\ast B.
\]
The representation of the reduced transfer function above is under the assumption that ${\mathcal V}$ and ${\mathcal W}$ are of equal dimension.

More generally, let us consider general $\cL_\infty$-functions in the framework of \eqref{eq:transfer_func}.
We define the reduced function by
\begin{equation}\label{eq:transfer_func_red}
	\widetilde{H} : \Omega \rightarrow {\mathbb C}^{p\times m},\quad	\widetilde{H}(s) := \widetilde{C}(s)  \widetilde{D}(s)^{-1} \widetilde{B}(s),
\end{equation}
where
\begin{subequations}\label{eq:middle_factor}
\begin{align}
	 \widetilde{B}(s)  &:=  f_1(s) \widetilde{B}_1 + \dots + f_{\kappa_B}(s) \widetilde{B}_{\kappa_B},& 
	 \quad {\widetilde B}_j &:= W^\ast B_j,& \quad &j = 1,\,\dots,\,\kappa_B, \\
	 \widetilde{C}(s)  &:=  g_1(s) \widetilde{C}_1 + \dots + g_{\kappa_C}(s) \widetilde{C}_{\kappa_C},& 
	 \quad {\widetilde C}_j &:= C_j V,& \quad &j = 1,\,\dots,\,\kappa_C, \\
	 \widetilde{D}(s)  &:=  h_1(s) \widetilde{D}_1 + \dots + h_{\kappa_D}(s) \widetilde{D}_{\kappa_D},&
	 \quad {\widetilde D}_j &:= W^\ast D_j V,& \quad &j = 1,\,\dots,\,\kappa_D.
\end{align}
\end{subequations}
Throughout the rest of this work, we focus on matrices {$V,\,W \in {\mathbb C}^{n\times \tilde{n}}$} with $\tilde{n} \ll n$ whose columns span the subspaces {${\mathcal V},\, {\mathcal W}$}, respectively. 
Furthermore, in what follows, we always assume that the subspaces ${\mathcal V},\, {\mathcal W}$ are such that
$\widetilde{H}$ is {well-defined and} bounded on the imaginary axis {and that} $\widetilde{D}$ is invertible {almost everywhere} on the imaginary axis.
The following result is fundamental to our approach. It is a special case of \cite[Theorem 1]{Gugercin2009}.

\begin{theorem}\label{thm:transfer_interpolate}
\label{thm1}
Let $\mu\in \mathbb{C}$ be such that $C(\mu)$, $D(\mu)$, and $B(\mu)$ are analytic and both 
$D(\mu)$and $\widetilde{D}(\mu)$ are invertible. Suppose also that $b\in \mathbb{C}^{m}$ and 
$c\in \mathbb{C}^{p}$ are given nonzero vectors. Then the following statements hold:
\begin{enumerate}
 	\item[\bf (i)] If $D(\mu)^{-1}B(\mu)b \in {\rm Col}(V)$, then $H(\mu)b = \widetilde{H}(\mu)b$;
	\item[\bf (ii)] If $\left(c^\ast C(\mu) D(\mu)^{-1} \right)^\ast \in {\rm Col}(W)$, then $c^\ast H(\mu) = c^\ast \widetilde{H}(\mu)$;
	\item[\bf (iii)] If $D(\mu)^{-1}B(\mu)b \in {\rm Col}(V)$ and $\left(c^\ast C(\mu) D(\mu)^{-1} \right)^\ast \in {\rm Col}(W)$,
	then $c^\ast H'\left(\mu\right)b = c^\ast \widetilde{H}'\left(\mu\right)b$.
\end{enumerate}
\end{theorem}

For the computation of the $\cL_\infty$-norm, we form subspaces ${\mathcal W},\, {\mathcal V}$ that give rise to
the Hermite interpolation of $\sigma(H(s))$ by $\sigma(\widetilde{H}(s))$ at some nodes $\mu_1,\, \dots,\, \mu_\ell$,
that is
\[
		\sigma(H(\mu_j))  =  \sigma\big(\widetilde{H}(\mu_j)\big)	\quad	{\rm and}	\quad	\sigma'(H(\mu_j))  =  \sigma'\big(\widetilde{H}(\mu_j)\big)
		\quad
		{\rm for } \quad  j = 1,\,\dots,\,\ell.
\] 
Theorem~\ref{thm:transfer_interpolate} above is helpful in this direction. It is immediate from part (i) of the theorem that if $D(\mu)^{-1}B(\mu)v \in {\mathcal V}$ for a right singular vector $v$ associated with $\sigma(H(\mu))$, then we have $\sigma(H(\mu)) \leq \sigma\big(\widetilde{H}(\mu)\big)$. The same conclusion can be drawn from part (ii) 
if $\left(w^\ast C(\mu) D(\mu)^{-1} \right)^\ast \in {\mathcal W}$ for a left singular vector $w$ associated 
with $\sigma(H(\mu))$.  Furthermore, it can be shown that 
if $D(\mu)^{-1}B(\mu)v_j \in {\mathcal V}$ for each right singular vector $v_j$ of $H(\mu)$ and
$\left(w^\ast_j C(\mu) D(\mu)^{-1} \right)^\ast \in {\mathcal W}$ for each left singular vector $w_j$ of $H(\mu)$,
then the equality $H(\mu) = \widetilde{H}(\mu)$ is attained (see Lemma~\ref{thm:sval_interpolate}, part (i) below),
implying $\sigma(H(\mu)) = \sigma\big(\widetilde{H}(\mu)\big)$. Additionally, 
$H(\mu) = \widetilde{H}(\mu)$ have the same right and left singular vectors $v,\, w$ corresponding to $\sigma(H(\mu)) = \sigma\big(\widetilde{H}(\mu)\big)$
and $D(\mu)^{-1}B(\mu)v \in {\mathcal V}, \left(w^\ast C(\mu) D(\mu)^{-1} \right)^\ast \in {\mathcal W}$. Consequently, part \textbf{(iii)}
of Theorem~\ref{thm:transfer_interpolate} leads to the desired Hermite interpolation property
\[
	\sigma'(H(\mu)) = \Real\left( w^\ast H'\left(\mu\right)v\right) = \Real\big( w^\ast \widetilde{H}'\left(\mu\right)v\big) = \sigma'\big(\widetilde{H}(\mu)\big),
\]
where the first and the third equality follow from the analytical formulas for the derivatives of singular value functions \cite{Lancaster1964, Bunse-Gerstner1991}.

An observation that enhances efficiency is that the singular vectors $v_j,\, w_j$ do not need to be calculated
explicitly. It is sufficient that we have 
\begin{align}
\label{eq:subspace_inc1} \left\{ D(\mu)^{-1}B(\mu)v_j \; \big| \; v_j \text{ is a right singular vector of } H(\mu) \right\}		&\subseteq	{\mathcal V}, \;\; {\rm and}	\\
\label{eq:subspace_inc2} \left\{\left(w^\ast_j C(\mu) D(\mu)^{-1} \right)^\ast \; \Big| \; w_j \text{ is a left singular vector of } H(\mu) \right\}		&\subseteq	{\mathcal W}
\end{align}
in order to obtain the Hermite interpolation property. Note that ${\mathcal V}$ and ${\mathcal W}$ must have the same
dimension, otherwise the middle factor $\widetilde{D}(s)$ of $\widetilde{H}(s)$ defined by \eqref{eq:middle_factor} is not square 
and Theorem~\ref{thm:transfer_interpolate} fails.
Clearly, the choices ${\mathcal V} = {\rm Col}(D(\mu)^{-1}B(\mu))$ and ${\mathcal W} = {\rm Col}\left( \left( C(\mu) D(\mu)^{-1} \right)^\ast \right)$
yield the desired inclusions \eqref{eq:subspace_inc1} and \eqref{eq:subspace_inc2}, but have different dimensions unless $m = p$.  
When $m < p$, we have
\begin{multline*}
  \left\{ \left(w^\ast_j C(\mu) D(\mu)^{-1} \right)^\ast \; \Big| \; w_j \text{ is a left singular vector of } H(\mu) \right\} 	\\
	 = \left\{ \left((H(\mu)v_j)^\ast C(\mu) D(\mu)^{-1} \right)^\ast \; \Big| \; v_j \text{ is a right singular vector of } H(\mu) \right\} \\	
	    	\subseteq  {\rm Col}\left(\left( C(\mu) D(\mu)^{-1} \right)^\ast H(\mu) \right),
\end{multline*}
so the subspaces
$
	{\mathcal V} = {\rm Col}(D(\mu)^{-1}B(\mu)) \quad	{\rm and}	\quad	{\mathcal W} = {\rm Col}\left( \left( C(\mu) D(\mu)^{-1} \right)^\ast H(\mu) \right) 
$
have equal dimension and satisfy \eqref{eq:subspace_inc1} and \eqref{eq:subspace_inc2}, respectively.
Similarly, when $m > p$, it can be deduced that the subspaces ${\mathcal V} = {\rm Col} \left( D(\mu)^{-1}B(\mu) H(\mu)^\ast \right)$ 
and ${\mathcal W} = {\rm Col}\left( \left( C(\mu) D(\mu)^{-1} \right)^\ast \right)$ are of equal dimension, and satisfy 
\eqref{eq:subspace_inc1} and \eqref{eq:subspace_inc2}.

The subspace method is described below in Algorithm \ref{alg1}. It generates matrices $V_r,\, W_r$ and acts on the subspaces 
${\rm Col}(V_r),\, {\rm Col}(W_r)$ of growing dimension as $r$ increases for $r = 1,\, 2,\, \dots$. 
In the description, the notation $\widetilde{H}_r(s)$ refers to the reduced function $\widetilde{H}(s)$ defined as in \eqref{eq:transfer_func_red} and \eqref{eq:middle_factor}, but with the particular choices $V = V_r$ and $W = W_r$. 
Thus, at iteration $r$ on line 10, the algorithm maximizes $\sigma\big( \widetilde{H}_{r-1}(s) \big)$ over the imaginary axis and retrieves the global maximizer $\ri \omega_r$. Then, it
expands the subspaces ${\rm Col}(V_{r-1}), {\rm Col}(W_{r-1})$ and thus forms $\widetilde{H}_r(s)$
such that the Hermite interpolation properties
$
	\sigma(H(\ri\omega_r))  =  \sigma\big(\widetilde{H}_r(\ri\omega_r)\big),\,\sigma'(H(\ri\omega_r))  =  \sigma'\big(\widetilde{H}_r(\ri\omega_r)\big)
$
hold. In practice we observe that Algorithm \ref{alg1} converges to a local maximizer of $\sigma(\omega)$ (that is
not necessarily a global maximizer) at a superlinear rate of convergence. The next section is devoted to a formal proof of this superlinear rate of convergence. Numerical experiments showing this convergence are reported in Section \ref{sec:experiment}.

\begin{algorithm}[tb]
 \begin{algorithmic}[1]
 
\REQUIRE{matrices $B_1,\,\dots,\,B_{\kappa_B} \in \mathbb{C}^{n\times m}$, $C_1,\,\dots,\,C_{\kappa_C} \in \mathbb{C}^{p\times n}$, $D_1,\,\dots,\,D_{\kappa_D} \in \mathbb{C}^{n\times n}$ and functions $f_1,\,\ldots,\,f_{\kappa_B},\,g_1,\,\ldots,\,g_{\kappa_C},h_1,\,\ldots,\,h_{\kappa_D}$ as in \eqref{eq:matsum}.}
\ENSURE{the $\cL_\infty$-norm of $H \in \cL_{\infty}^{p \times m}$ with $H$ as in \eqref{eq:transfer_func} and \eqref{eq:matsum}.}
\STATE $\omega_{1} \gets$ a random number in ${\mathbb R}.$
\IF{$m = p$}
	\STATE $V_1 \gets  D(\ri\omega_1 )^{-1}B(\ri\omega_1)\quad \text{and}\quad W_1 \gets \left( C(\ri\omega_1)D(\ri \omega_1 )^{-1} \right)^\ast$.
\ELSIF{$m < p$}
	\STATE $V_1 \gets  D(\ri\omega_1 )^{-1}B(\ri\omega_1)\quad \text{and}\quad W_1 \gets \left( C(\ri\omega_1)D(\ri \omega_1 )^{-1} \right)^\ast H(\ri \omega_1)$.
\ELSE
	\STATE $V_1 \gets  D(\ri\omega_1 )^{-1}B(\ri\omega_1) H(\ri \omega_1)^\ast \quad \text{and}\quad W_1 \gets \left( C(\ri\omega_1)D(\ri \omega_1 )^{-1} \right)^\ast$.
\ENDIF

	\FOR{$r = 2,\,3,\,\dots$}
	
	\STATE  Form $\widetilde{H}_{r-1}$ as in \eqref{eq:transfer_func_red} and \eqref{eq:middle_factor} and set $\omega_r \gets \arg\max_{\omega \in {\mathbb R} \cup \{\infty\}} \sigma(\widetilde{H}_{r-1} (\ri \omega))$.
	\IF{$m = p$}
		\STATE $\widetilde{V}_r \gets  D(\ri\omega_r )^{-1}B(\ri\omega_r)\quad\text{and}\quad \widetilde{W}_r \gets \left( C(\ri\omega_r)D(\ri \omega_r )^{-1} \right)^\ast$.
	\ELSIF{$m < p$}
		\STATE $\widetilde{V}_r \gets  D(\ri\omega_r )^{-1}B(\ri\omega_r)\quad\text{and}\quad \widetilde{W}_r \gets \left( C(\ri\omega_r)D(\ri \omega_r )^{-1} \right)^\ast H(\ri \omega_r)$.
	\ELSE
		\STATE $\widetilde{V}_r \gets  D(\ri\omega_r )^{-1}B(\ri\omega_r) H(\ri \omega_r)^\ast \quad\text{and}\quad \widetilde{W}_r \gets \left( C(\ri\omega_r)D(\ri \omega_r )^{-1} \right)^\ast$.
	\ENDIF
	\STATE $V_r \gets \operatorname{orth}\left(\begin{bmatrix} V_{r-1} & \widetilde{V}_r \end{bmatrix}\right)
		\quad \text{and}\quad W_r \gets \operatorname{orth}\left(\begin{bmatrix} W_{r-1} & \widetilde{W}_r \end{bmatrix}\right).$
\ENDFOR
 \end{algorithmic}
\caption{Subspace method for the computation of the $\cL_\infty$-norm}
\label{alg1}
\end{algorithm}

\begin{remark}
The procedure described in Algorithm~\ref{alg1} resembles the reduced
basis approach for model order reduction of parametrized systems, see,
e.\,g., \cite{Haasdonk16}. The key ingredients are projection onto a
subspace, solving the resulting low-dimensional problem, a
subprocedure to maximize (minimize) a desired quantity for the
reduced parametrized system, and expanding the subspace by a snapshot
of the full-order problem at the argmax/min returned by the subprocedure. 
As all these ingredients are used in Algorithm~\ref{alg1}, it can be considered as a
reduced basis method. 
\end{remark}

\begin{remark}
In Algorithm~\ref{alg1} the subspaces from all of the previous iterations are kept.
An alternative would be to keep the subspaces only from the last two iterations. 
The rate of convergence analysis in the next section also applies to this variant,
since that analysis (specifically Theorem \ref{thm:quad_conv})
relies on the interpolation properties only at the last two iterates.
Thus, the variant with only subspaces from the last two iterations is also guaranteed to
converge at a superlinear rate, which we observe in practice. {However, for the numerical 
experiments discussed in Section \ref{sec:experiment}, all the previous subspaces are kept.}
This results in better global convergence properties, and usually avoidance of stagnation at a local maximizer of $\sigma(\cdot)$
that is not a global maximizer. The cost of keeping additional subspaces is usually small,
{because the algorithm often needs fewer iterations than 
the variant that uses only the subspaces from the last two iterations. Typically it converges quickly up to prescribed tolerances in less} 
than 10 iterations (see the numerical results in Sections \ref{sec:nexam_desc} and \ref{sec:nexam_delay}).
\end{remark}

\section{Rate of Convergence Analysis}\label{sec:convergence}
In this section, we prove that the aforementioned Hermite interpolation properties of the subspace 
method lead to a superlinear convergence with respect to the subspace dimension, under the
assumption that the method converges locally. The argument revolves around the singular value 
functions $\sigma(\omega) := \sigma(H(\ri\omega))$ and $\sigma_r(\omega) := \sigma\big(\widetilde{H}_r(\ri\omega)\big)$. 
Occasionally, the second largest singular values of $H(\ri\omega)$ and $\widetilde{H}_r(\ri\omega)$
are also referred, which we denote by $\underline{\sigma}(\omega)$ and $\underline{\sigma}_r(\omega)$,
respectively. 
When $\min \{ m, p \} = 1$, {then} we define $\underline{\sigma}(\omega) = \underline{\sigma}_r(\omega) = 0$ for all $\omega \in {\mathbb R}$.
We first formally state and prove Hermite interpolation properties of the singular value functions. 

\begin{lemma}\label{thm:sval_interpolate}
The following statements hold regarding Algorithm~\ref{alg1} for $k = 1,\,\ldots,\,r$:
\begin{enumerate}
	\item[\bf (i)] $H(\ri\omega_k) = \widetilde{H}_r(\ri\omega_k)$;
	\item[\bf (ii)] $\sigma(\omega_k) = \sigma_r(\omega_k)$ and $\underline{\sigma}(\omega_k) = \underline{\sigma}_r(\omega_k)$;
	\item[\bf (iii)] If $\sigma(\omega_k)$ is simple, then $\sigma(\omega),\,\sigma_r(\omega)$ are differentiable
	at $\omega_k$ and $\sigma'(\omega_k) = \sigma'_r(\omega_k)$.
\end{enumerate}
\end{lemma}

\begin{proof}
\begin{enumerate}
\item[\textbf{(i)}] When $m \leq p$, for each $k \in \{1,\,\dots,\,r\}$, we have $\Col(D(\ri\omega_k)^{-1} B(\ri\omega_k)) \subseteq \Col(V_r)$,
due to lines 3, 5, 12, 14, and 18 of Algorithm~\ref{alg1}.
Thus, $D(\ri\omega_k)^{-1} B(\ri\omega_k) e_j \in \Col(V_r)$ for $j = 1,\,\dots,\,m$. It follows that $H(\ri\omega_k) e_j = \widetilde{H}_r(\ri\omega_k) e_j$ from part (i) of Theorem~\ref{thm:transfer_interpolate} for $j = 1,\,\dots,\,m$, that is $H(\ri\omega_k) = \widetilde{H}_r(\ri\omega_k)$. 
On the other hand, when $m > p$, for each $k \in \{1,\,\dots,\,r\}$, the inclusion $\Col \left( \left( C(\ri \omega_k)D(\ri \omega_k )^{-1} \right)^\ast \right) \subseteq \Col(W_r)$ follows from lines 7, 16, and 18 of Algorithm~\ref{alg1}. Consequently, $\left( e_j^\ast C(\ri \omega)D(\ri \omega_k)^{-1} \right)^\ast  \in \Col(W_r)$,
so $e_j^\ast H(\ri\omega_k) = e_j^\ast \widetilde{H}_r(\ri\omega_k)$ by part (ii) of Theorem~\ref{thm:transfer_interpolate} for each $j = 1,\,\dots,\, p$, that is $H(\ri\omega_k) = \widetilde{H}_r(\ri\omega_k)$. 
\item[\textbf{(ii)}] This is immediate from part (i). 
\item[\textbf{(iii)}] Suppose that $\sigma(\omega_k)$ is simple for a particular $k \in \{1,\,\dots,\,r\}$. This implies that
$\sigma(\omega),\,\sigma_r(\omega)$ are differentiable at $\omega_k$ \cite{Rellich1969}. 
The left and right singular vectors corresponding to
$\sigma(\omega_k)$ and $\sigma_r(\omega_k)$ are the same, since
$H(\ri\omega_k) = \widetilde{H}_r(\ri\omega_k)$ due to part (i). Let
us denote them by $w\in {\mathbb C}^p$ and $v \in {\mathbb C}^m$,
respectively, and w.\,l.\,o.\,g., assume these are unit vectors. Suppose $m \leq p$. In this case, $\Col( D(\ri\omega_k)^{-1} B(\ri\omega_k)) \subseteq \Col(V_r)$ and $\Col\left(\left( C(\ri\omega_k)D(\ri \omega_k )^{-1} \right)^\ast H( \ri \omega_k ) \right) \subseteq \Col(W_r)$, so we have $D(\ri\omega_k)^{-1} B(\ri\omega_k) v \in \Col(V_r)$ and 
\[
	\left( C(\ri\omega_k)D(\ri \omega_k )^{-1} \right)^\ast H( \ri \omega_k )  v  = \sigma(\omega_k) \left( w^\ast C(\ri\omega_k)D(\ri \omega_k )^{-1} \right)^\ast \in \Col(W_r).
\]
When we have $m > p$, the inclusions $\Col( D(\ri\omega_k)^{-1} B(\ri\omega_k) H(\ri\omega_k)^\ast ) \subseteq \Col(V_r)$ 
and $\Col\left(\left( C(\ri\omega_k) D(\ri \omega_k )^{-1} \right)^\ast \right) \subseteq \Col(W_r)$ hold. This implies
that $\left( w^\ast C(\ri\omega_k)D(\ri \omega_k )^{-1} \right)^\ast \in \Col(W_r)$ and
\[
	D(\ri\omega_k)^{-1} B(\ri\omega_k) H(\ri\omega_k)^\ast w = \sigma(\omega_k) D(\ri\omega_k)^{-1} B(\ri\omega_k) v \in \Col(V_r).
\]
In both cases, part (iii) of Theorem~\ref{thm:transfer_interpolate} yields $w^\ast H'(\ri \omega_k) v = w^\ast \widetilde{H}'_r(\ri \omega_k) v$.
Finally, by exploiting the analytical formulas for the derivatives of singular value functions \cite{Lancaster1964, Bunse-Gerstner1991}, we deduce
\[
	\sigma'(\omega_k) = \Real( w^\ast H'(\ri\omega_k) v ) = \Real\big( w^\ast \widetilde{H}'_r(\ri\omega_k) v \big) = \sigma'_r(\omega_k).
\]
\end{enumerate}
\end{proof}


The next result concerns how accurately $\sigma''_r(\cdot)$ approximates $\sigma''(\cdot)$ 
at $\omega_r$. We view $\omega_k$ for every $k > 1$ as a function of the initial point $\omega_1$ for the next
result and the subsequent rate of convergence result. A consequence is that the function $\sigma_k(\cdot)$
also depends on $\omega_1$. Furthermore, in what follows, for a given bounded interval ${\mathcal I} \subseteq {\mathbb R}$, 
we consider $\omega_1$ such that $\omega_k \in {\mathcal I}$ for each {$k \ge 1$}. Due to the analyticity 
of the function $H$ on the imaginary axis (recall that $H \in  \cL_\infty^{p \times m}$), there exists a Lipschitz constant $\eta_{\mathcal I} > 0$ such that
\begin{equation}\label{eq:Lipschitz_cond}
		\| H(\ri \omega) - H(\ri \widetilde{\omega}) \|_2	\leq	\eta_{\mathcal I} |\omega - \widetilde{\omega} |
		\quad\quad
		\forall \omega,\, \widetilde{\omega} \in {\mathcal I}.
\end{equation}
Additionally, for a given $\eta \geq \eta_{\mathcal I}$ and $\varphi > 0$, we consider $\omega_1$ such that
\begin{equation}\label{eq:Lipschitz_cond2}
	\| \widetilde{H}^{(j)}_k(\ri \omega) - \widetilde{H}^{(j)}_k(\ri \widetilde{\omega}) \|_2	\leq	\eta |\omega - \widetilde{\omega} |	\quad\quad
		\forall \omega,\, \widetilde{\omega} \in {\mathcal I}
\end{equation}
for $k \geq 1$ and $j = 0, 1, 2$ as well as
\begin{equation}\label{eq:bounded_derivatives}
	\| \widetilde{H}^{(j)}_k (\ri \omega) \|_2 \leq \varphi	\quad\quad	\forall \omega \in {\mathcal I}
\end{equation}
for $k \geq 1$ and $j =  1, 2$. Condition (\ref{eq:Lipschitz_cond}), in particular the existence of the constant $\eta_{\mathcal I}$,  
is a simple consequence of the analyticity of $H(\ri \cdot)$ and the boundedness of ${\mathcal I}$, whereas conditions
(\ref{eq:Lipschitz_cond2}) and (\ref{eq:bounded_derivatives}) are assumptions, which are typically satisfied in practice
because of the interpolation properties between $H(\ri \cdot)$ and $\widetilde{H}_k(\ri \cdot)$. These conditions
imply the Lipschitz continuity of $\sigma''(\cdot)$ and $\sigma''_k(\cdot)$ on ${\mathcal I}$ with Lipschitz constants 
independent of $\omega_1$, which is established and exploited by the proof of the next lemma.

\begin{lemma}\label{thm:accuracy_sderivative}
For a given $\zeta \in {\mathbb R}^+$, an integer $r \geq 2$, a bounded interval ${\mathcal I} \subseteq {\mathbb R}$, $\varphi \in {\mathbb R}^+$, 
and $\eta \geq \eta_{\mathcal I}$, where $\eta_{\mathcal I}$ is as in \eqref{eq:Lipschitz_cond}, suppose $\omega_1$ 
is chosen in a way so that $\omega_k \in {\mathcal I}$ for {$k \ge 1$}, conditions \eqref{eq:Lipschitz_cond2} 
and \eqref{eq:bounded_derivatives} hold, as well as
\begin{equation}\label{eq:separation}
	\sigma(\omega_r) - \underline{\sigma}(\omega_r) \geq \zeta  \geq  c \eta | \omega_r  -  \omega_{r-1} |
\end{equation}
for some constant $c > 2$. Then we have
\[
	|\sigma''(\omega_r) - \sigma''_r(\omega_r)| 
		 \leq  \mu  | \omega_r  -  \omega_{r-1} |	
\]
for some constant $\mu$ independent of $\omega_1$.
%
%
%
\end{lemma}

\begin{proof}
We start by establishing the simplicity of $\sigma(\omega)$ and $\sigma_r(\omega)$ on the closed interval with end-points 
$\omega_{r-1},\, \omega_r$, which we denote with ${\mathcal I}_r$. To this end, for each $\omega \in {\mathcal I}_r$, we have
\begin{align*}
| \sigma(\omega_r) - \sigma(\omega) | \leq \| H(\ri \omega_r) - H(\ri \omega) \|_2	\leq	\eta | \omega_r - \omega |	  \leq &	\eta | \omega_r - \omega_{r-1} |	 \quad \text{and}	\\
| \underline{\sigma}(\omega_r) - \underline{\sigma}(\omega) | \leq \| H(\ri \omega_r) - H(\ri \omega) \|_2	\leq	\eta | \omega_r - \omega |	\leq &	\eta | \omega_r - \omega_{r-1} |,
\end{align*}
due to Weyl's theorem \cite[Theorem 4.3.1]{Horn1985} and inequality \eqref{eq:Lipschitz_cond} regarding the Lipschitz continuity of $H$.
Hence we have
\[
		\sigma(\omega)	\geq		\sigma(\omega_r)	-	\eta | \omega_r - \omega_{r-1} |			 \quad 	\text{and}		\quad
		\underline{\sigma}(\omega)	\leq		\underline{\sigma}(\omega_r)		+	\eta | \omega_r - \omega_{r-1} |	,	
\]
that is
\[
\sigma(\omega)	-	\underline{\sigma}(\omega)		\geq	\left\{ \sigma(\omega_r)	-	\underline{\sigma}(\omega_r) \right\} 	-	2\eta | \omega_r - \omega_{r-1} |
												\geq	 (c-2)	 \eta | \omega_r - \omega_{r-1} |	> 0.
\]
Above, the second inequality follows from \eqref{eq:separation}. This shows that $\sigma(\omega)$ is simple for each $\omega \in {\mathcal I}_r$.
Furthermore, by part (ii) of Lemma~\ref{thm:sval_interpolate}, we have $\sigma_r(\omega_r) = \sigma(\omega_r)$ and 
$\underline{\sigma}_r(\omega_r) = \underline{\sigma}(\omega_r)$. An analogous argument with $\sigma_r(\cdot)$ taking the
role of $\sigma(\cdot)$ also shows the simplicity of $\sigma_r(\omega)$ for each $\omega \in {\mathcal I}_r$.
It follows that both $\sigma(\cdot)$ and $\sigma_r(\cdot)$ are analytic on ${\mathcal I}_r$.

To relate the second derivatives, we exploit part (iii) of Lemma~\ref{thm:sval_interpolate}, in particular
$
\sigma'(\omega_k)  =  \sigma'_r(\omega_k)  
$
for $k = r-1,\, r$. These interpolation properties imply
\[
\sigma'' (\xi ) \left( \omega_r- \omega_{r-1} \right)   =  \sigma'(\omega_r) - \sigma'(\omega_{r-1}) 
									     =  \sigma'_r(\omega_r) -  \sigma'_r(\omega_{r-1}) 
									     =  \sigma''_r \big(\widehat{\xi} \big) \left( \omega_r - \omega_{r-1} \right)
\]
for some $\xi,\,\widehat{\xi} \in {\mathcal I}_r$ leading to
\begin{equation}\label{zero}
	 \sigma'' (\xi)  -  \sigma''_r (\widehat{\xi})  = 0.
\end{equation}
Moreover, the second derivatives of $\sigma(\cdot),\, \sigma_r(\cdot)$ are Lipschitz continuous in ${\mathcal I}_r$, so there 
exist positive constants $\gamma_1, \gamma_2$ such that
\begin{equation}\label{eq:Lipschitz_sder}
\begin{split}
\big| \sigma''_r  (\omega_r) -  \sigma''_r  \big(\widehat{\xi}\big) \big|  &\leq \gamma_1 \big| \widehat{\xi} -\omega_r \big|   \leq  \gamma_1  | \omega_r - \omega_{r-1} |  \quad \text{and} \\
| \sigma'' (\omega_r) - \sigma'' (\xi) | &\leq \gamma_2 | \xi -\omega_r |   \leq  \gamma_2 | \omega_r - \omega_{r-1} |.
\end{split}
\end{equation}
We claim that the Lipschitz constant $\gamma_1$ can be expressed solely in terms of $\eta,\, \varphi,\, \zeta$
(satisfying (\ref{eq:Lipschitz_cond2}), (\ref{eq:bounded_derivatives}), (\ref{eq:separation}), respectively).

To see this, let us denote a unit eigenvector corresponding to the largest eigenvalue of
\[
						\left[
							\begin{array}{cc}
								0							&	\widetilde{H}_r(\ri \omega)	\\
								\left[ \widetilde{H}_r(\ri \omega) \right]^\ast	&	0
							\end{array}
						\right]
\] 
by $v_r(\omega)$, and a unit eigenvector corresponding to the $j$-th largest eigenvalue $\lambda_{r,j}(\omega)$ 
of this matrix by $v_{r,j}(\omega)$. Then the claim is evident from the analytical expression \cite{Lan64}
\begin{align*}
	\sigma''_r  (\omega)	= & \; 	v_r(\omega)^\ast	
						\left[
							\begin{array}{cc}
								0							&	\widetilde{H}''_r(\ri \omega)	\\
								\left[ \widetilde{H}''_r(\ri \omega) \right]^\ast	&	0
							\end{array}
						\right]
						v_r(\omega)
								\\
						+ & \; 	2
						\sum_{j=2}^{2r}	\frac{1}{\sigma_r(\omega) - \lambda_{r,j}(\omega)}	
								\left| v_{r,j}(\omega)^\ast
										\left[
											\begin{array}{cc}
												0							&	\widetilde{H}'_r(\ri \omega)	\\
												\left[ \widetilde{H}'_r(\ri \omega) \right]^\ast	&	0
											\end{array}
										\right]
									v_r(\omega)	\right|^2,
\end{align*}
where $\sigma_r(\cdot),\, \lambda_{r,j}(\cdot),\, v_r(\cdot),\, v_{r,j}(\cdot),\, \widetilde{H}'_r(\ri \cdot)$, and $\widetilde{H}''_r(\ri \cdot)$ 
are Lipschitz continuous on $\mathcal{I}_r$ with Lip\-schitz constants depending on $\eta$ only.
Here we remark that the terms $\sigma_r(\omega) - \lambda_{r,j}(\omega)$ can be bounded from below
by a quantity solely depending on $\zeta$, because of the interpolation properties $\sigma(\omega_r)  =  \sigma_r(\omega_r)$,
$\underline{\sigma}(\omega_r)  =  \underline{\sigma}_r(\omega_r)$ and assumption (\ref{eq:separation}).
Furthermore we exploit the fact that if $f,\,g$ are Lipschitz continuous functions with Lipschitz constants $\beta_1,\, \beta_2$
on a closed interval, then $fg$ is also Lipschitz continuous with Lipschitz constant $\beta_1 g_\ast + \beta_2 f_\ast$,
where $f_\ast$, $g_\ast$ are the maximum values of $f$, $g$ attained on the interval.
Similarly, the Lipschitz constant $\gamma_2$ in (\ref{eq:Lipschitz_sder}) can be expressed in terms 
of $\eta,\,\zeta$ and an upper bound on $\| H^{(j)}(\ri \omega) \|_2$ for $j = 1,\,2$ and for all $\omega \in {\mathcal I}$.
Finally, equation \eqref{zero} and inequalities \eqref{eq:Lipschitz_sder} yield
\begin{align*}
|  \sigma'' (\omega_r)  -  \sigma''_r (\omega_r)  | &=   |  \sigma''(\omega_r)   -  \sigma''_r(\omega_r) +  \sigma''_r (\widehat{\xi})  -   \sigma'' (\xi) | \\
		& \leq | \sigma'' (\omega_r) -  \sigma'' (\xi) | + \big| \sigma''_r (\omega_r ) - \sigma''_r \big(\widehat{\xi}\big) \big|  \leq (\gamma_1 + \gamma_2) | \omega_r - \omega_{r-1} |,
\end{align*}
hence the result follows.
\end{proof}

The main result presented next assumes $\omega_{r-1},\, \omega_r,\, \omega_{r+1}$ are sufficiently close to a local 
maximizer of $\sigma(\cdot)$ for certain values of $\omega_1$ and a given $r \in \N$. This is a 
convergence assumption which we observe in practice.

\begin{theorem}[Local superlinear convergence]\label{thm:quad_conv}
Let $\omega_\ast$ be a local maximizer of $\sigma\left(\omega\right)$ such that $\sigma ( \omega_\ast )$ is simple,
and $\sigma''(\omega_\ast) \neq 0$. Furthermore, let $\zeta := \sigma(\omega_\ast) - \underline{\sigma}(\omega_\ast)$. 
For a given integer $r \geq 2$, a bounded interval ${\mathcal I} \subseteq {\mathbb R}$ containing $\omega_\ast$ {in its interior}, 
$\varphi \in {\mathbb R}^+$, and $\eta \geq \eta_{\mathcal I}$, where $\eta_{\mathcal I}$ is as in \eqref{eq:Lipschitz_cond}, 
suppose $\omega_1$ is chosen in a way so that $\omega_k \in {\mathcal I}$ for {$k \ge 1$},
conditions \eqref{eq:Lipschitz_cond2}, \eqref{eq:bounded_derivatives} hold, 
and $\delta := \max\left\{ | \omega_{r+1}  -  \omega_{\ast} |, | \omega_{r}  -  \omega_{\ast} |, | \omega_{r-1}  -  \omega_{\ast} | \right\}$
is sufficiently small, in particular 
\begin{equation}\label{eq:separation2}
		 \zeta	\geq   8 \eta \delta.
\end{equation}
%
Then we have
\[
	\frac{  | \omega_{r+1} - \omega_\ast |  }{  | \omega_r - \omega_\ast | \cdot {\max} \left\{  | \omega_{r-1} - \omega_\ast | , | \omega_r - \omega_\ast | \right\} }
		 \leq  \nu		
\]
for some constant $\nu$ independent of $\omega_1$.
\end{theorem}

\begin{proof}
The proof is split into two parts. In the first part, we deduce the analyticity of the singular value 
functions $\sigma(\cdot)$, $\sigma_r(\cdot)$ on ${\mathcal I}(\omega_\ast,\delta) := [\omega_\ast-\delta, \omega_\ast+\delta]$, 
bound their second derivatives from below and the third derivative of
$\sigma_r(\cdot)$ from above uniformly on ${\mathcal I}(\omega_\ast,\delta)$
by quantities that do not depend on $\omega_1$. Then the second part makes use of these uniform bounds to 
relate $| \omega_{r+1}  -  \omega_{\ast} |$ with $| \omega_{r}  -  \omega_{\ast} |, | \omega_{r-1}  -  \omega_{\ast} |$
and {conclude} a superlinear rate of convergence.

\emph{Part 1:} We first show the analyticity of $\sigma(\cdot)$ and $\sigma_r(\cdot)$ on ${\mathcal I}(\omega_\ast,\delta)$.
Condition \eqref{eq:separation2} together with Weyl's theorem \cite[Theorem 4.3.1]{Horn1985} ensures that 
\[
	\sigma(\omega)	-	\underline{\sigma}(\omega)	\geq		\left\{ \sigma(\omega_\ast)	-	\underline{\sigma}(\omega_\ast) \right\} - 2\eta \delta
												\geq		3\zeta/4
\]
for each $\omega \in {\mathcal I}(\omega_\ast,\delta)$, meaning $\sigma(\omega)$ is  simple on this interval. Moreover,
\[
	\sigma_r(\omega_r) - \underline{\sigma}_r(\omega_r)  = \sigma(\omega_r) - \underline{\sigma}(\omega_r) \geq 	3 \zeta/4.
\]
But $|{\mathcal I}(\omega_\ast,\delta)| = 2\delta$, so, by Weyl's theorem, we also have
\[
	\sigma_r(\omega) - \underline{\sigma}_r(\omega) \geq \left\{ \sigma_r(\omega_r) - \underline{\sigma}_r(\omega_r) \right\}	-	4\eta\delta
											\geq	\zeta/4
\] 
for all $\omega \in {\mathcal I}(\omega_\ast,\delta)$.
Consequently, $\sigma(\cdot)$ and $\sigma_r(\cdot)$ are analytic on ${\mathcal I}(\omega_\ast,\delta)$.

Secondly, we show that the second derivatives of $\sigma(\cdot)$ and $\sigma_r(\cdot)$ are
bounded away from zero on ${\mathcal I}(\omega_\ast,\delta)$. For the
former, w.\,l.\,o.\,g.\, due  
to $\sigma''(\omega_\ast) \neq  0$, we simply consider $\delta$ small enough (much smaller than $|\sigma''(\omega_\ast) |$) so that
\[
	 |\sigma''(\omega)| \geq  \ell_1	\quad		\forall \omega \in {\mathcal I}(\omega_\ast,\delta)
\]
for some constant $\ell_1 \gg \delta > 0$. For the latter,
\[
		3 \eta | \omega_r - \omega_{r-1} |  \leq  3 \eta ( | \omega_r - \omega_\ast | + |\omega_\ast - \omega_{r-1} |  )
				 \leq  6\eta \delta \leq \frac{3\zeta}{4}  \leq  \sigma(\omega_r) - \underline{\sigma}(\omega_r),
\]
so Lemma~\ref{thm:accuracy_sderivative} implies
\begin{equation}\label{eq:sder_proximity}
	| \sigma''_r(\omega_r) - \sigma''(\omega_r) |	 	\leq 		\mu	| \omega_r  -  \omega_{r-1}  |	 \leq  
			2\mu  \max \left\{ | \omega_r - \omega_\ast |, | \omega_{r-1} - \omega_\ast | \right\}  \leq  2\delta \mu.
\end{equation}
That is, $\: | \sigma_r''(\omega_r) | \geq |\sigma''(\omega_r)| - 2\delta \mu \geq \ell_1 - 2\delta \mu$. It follows that
\[
	 |\sigma_r''(\omega)|  \geq  \ell_2	\quad		\forall \omega \in {\mathcal I}(\omega_\ast,\delta)
\]
for another constant $\ell_2 > 0$. 

Thirdly, we show that {the absolute value of the third derivative} of $\sigma_r(\cdot)$ is bounded from above
by a quantity that does not depend on $\omega_1$ on ${\mathcal I}(\omega_\ast,\delta)$.
Repeating the arguments in the proof of Lemma \ref{thm:accuracy_sderivative},
the second derivative $\sigma''_r(\cdot)$ is Lipschitz continuous on ${\mathcal I}(\omega_\ast,\delta)$
with a Lipschitz constant $\gamma_1$ that depends on $\eta$, $\varphi$ and $\zeta$ only. This in turn, together
with the analyticity of $\sigma_r(\cdot)$ on ${\mathcal I}(\omega_\ast,\delta)$, {implies}
\[
	|\sigma_r'''(\omega)|  \leq  \gamma_1	\quad		\forall \omega \in {\mathcal I}(\omega_\ast,\delta).
\]

\emph{Part 2:} We express $| \omega_{r+1}  -  \omega_\ast |$ in terms of $| \omega_r  -  \omega_\ast |$ and $| \omega_{r-1} - \omega_\ast |$, and conclude with the superlinear convergence result as desired. Analyticity of $\sigma(\omega)$ implies
\[
	0  = \sigma'(\omega_\ast) =  \sigma'(\omega_r) +  \int_0^1 \sigma'' \left(\omega_r + t (\omega_\ast - \omega_r )\right) \left( \omega_\ast - \omega_r \right) \mathrm{d}t.
\]
In the last equation, we employ $\sigma'(\omega_r) = \sigma'_r(\omega_r)$ (part (iii) of Lemma~\ref{thm:sval_interpolate}),
divide both sides by $\sigma''(\omega_r)$, and reorganize to obtain
\begin{equation}\label{eq:taylor_organized}
	0 =  \frac{\sigma_r'(\omega_r)}{\sigma''(\omega_r)} + \left(\omega_\ast - \omega_r \right) +
	\frac{1}{\sigma''(\omega_r)}
			\int_0^1 \left[	\sigma'' \left(\omega_r + t (\omega_\ast - \omega_r )\right) - \sigma''(\omega_r) 	\right]
							\left(\omega_\ast - \omega_r \right) \mathrm{d}t.
\end{equation}
In what follows, we exploit $\sigma''(\omega_r) \approx \sigma''_r(\omega_r)$ as a consequence of Lemma~\ref{thm:accuracy_sderivative},
and  $\sigma'_r(\omega_r) / \sigma''_r(\omega_r) \approx -(\omega_{r+1} - \omega_r)$ as a consequence of $\sigma'_r(\omega_{r+1}) = 0$.
These observations lead us to $\sigma'_r(\omega_r) / \sigma''(\omega_r) + (\omega_\ast - \omega_r) \approx (\omega_\ast - \omega_{r+1})$
in (\ref{eq:taylor_organized}). Formally, an application of Taylor's theorem with Lagrange remainder to $\sigma_r'(\cdot)$ 
and optimality of $\omega_{r+1}$ with respect to $\sigma_r(\cdot)$ give rise to
\begin{equation*}
				0  =  \sigma'_r(\omega_{r+1}) =  \sigma'_r(\omega_r)  + \sigma''_r(\omega_r) (\omega_{r+1} - \omega_r)	
													+	\frac{\sigma'''_r(\xi)}{2} (\omega_{r+1} - \omega_r)^2,
\end{equation*}													
which can be rearranged as
\begin{equation}
\label{eq:relate_der}					
	\frac{ \sigma'_r(\omega_r) }{ \sigma''_r(\omega_r) }  =  
			- \left( \omega_{r+1} - \omega_r \right) -	\frac{\sigma'''_r(\xi)}{2 \sigma''_r(\omega_r) } (\omega_{r+1} - \omega_r)^2
\end{equation}
for some $\xi \in {\mathcal I}(\omega_\ast,\delta)$. By combining \eqref{eq:taylor_organized} and \eqref{eq:relate_der}, we deduce
\begin{multline*}
	 0 =  \left( \omega_\ast - \omega_{r+1} \right) + \left(\frac{1}{\sigma''(\omega_r)} - \frac{1}{\sigma''_r(\omega_r)} \right) \sigma'(\omega_r)
 		-	\frac{\sigma'''_r(\xi)}{2 \sigma''_r(\omega_r) } (\omega_{r+1} - \omega_r)^2	 \\
 		+ \frac{1}{\sigma''(\omega_r)} 
			\int_0^1 \left[	\sigma'' \left(\omega_r + t (\omega_\ast - \omega_r )\right) - \sigma''(\omega_r)	\right]	(\omega_\ast - \omega_r)		\mathrm{d}t,
 \end{multline*}  
implying
\begin{multline}\label{eq:taylor_simplified}
 	| \omega_\ast - \omega_{r+1}  |   \leq  
		\left| \frac{\sigma''_r(\omega_r) - \sigma''(\omega_r)}{\sigma''(\omega_r)\sigma''_r(\omega_r)} \right|  \left| \sigma'(\omega_r) \right|
 		+	\left| 	\frac{\sigma'''_r(\xi)}{2 \sigma''_r(\omega_r) } \right|  | \omega_{r+1} - \omega_r |^2 		\\	
		+	 	 \left|  \frac{\gamma}{2 \sigma''(\omega_r)}  \right| |\omega_\ast - \omega_r|^2
\end{multline}
where we used the Lipschitz continuity of $\sigma''(\cdot)$ on the interval ${\mathcal I}(\omega_\ast,\delta)$, in particular
we used the existence of a Lipschitz constant $\gamma > 0$ such that
\[
	| \sigma'' \left(\omega_r + t (\omega_\ast - \omega_r )\right) - \sigma''(\omega_r) |	\leq 	\gamma t |\omega_\ast - \omega_r |	\quad	\forall t \in [0,1].
\]
 
Finally, by Young's inequality, we have
$| \omega_{r+1}  -  \omega_r |^2 \leq 2 | \omega_{r+1} - \omega_\ast |^2 + 2| \omega_{r} - \omega_\ast |^2$. Thus, 
the expression (\ref{eq:taylor_simplified}) can be rewritten as
\begin{equation}\label{eq:rate_of_convergence_close}
 	(1 - c_1 \delta) | \omega_{r+1} - \omega_\ast |	 \leq		
	c_2 | \sigma''_r(\omega_r) - \sigma''(\omega_r) | | \sigma'(\omega_r) | +  c_3 | \omega_r - \omega_\ast|^2
\end{equation}
where $c_1 = \gamma_1/\ell_2$, $c_2 = 1/(\ell_1 \cdot \ell_2)$ and $c_3 = \gamma/(2\ell_1) + \gamma_1/\ell_2$.
The term $| \sigma''_r(\omega_r) - \sigma''(\omega_r) |$ on the right in (\ref{eq:rate_of_convergence_close}) is bounded from above by 
 ${\rm max} \left\{ | \omega_r - \omega_\ast |, | \omega_{r-1} - \omega_\ast | \right\}$ up to a constant by \eqref{eq:sder_proximity},
whereas the term $| \sigma'(\omega_r) |$ on the right is bounded from above by $| \omega_r - \omega_\ast|$ up to a constant by the 
mean value theorem. 
If $\delta$ is sufficiently small, the term on the left-hand side of \eqref{eq:rate_of_convergence_close} can be bounded from below by 
$c_4 | \omega_{r+1} - \omega_\ast |$ for some constant $c_4 > 0$. Hence, the result follows.
\end{proof}

\section{Numerical Experiments}\label{sec:experiment}
In this section, we report on the numerical results obtained by our MATLAB implementation of Algorithm~\ref{alg1}\footnote[1]{available from \url{http://www.math.tu-berlin.de/index.php?id=186267&L=1}}. We first describe a few important implementation details and the test setup. After that we report on the respective numerical results. 
\subsection{Implementation Details and Test Setup}
At each iteration of Algorithm~\ref{alg1}, the $\cL_\infty$-norm of a reduced function needs to be computed in line 10. This global nonconvex optimization problem is solved by means of the approach due to Boyd and Balakrishnan for transfer functions of linear state-space systems \cite{Boyd1990} (and \cite{BenSV12} for the case of descriptor systems), and by means of the algorithm in \cite{Mengi2014} for general $\cL_\infty$-functions. The Boyd-Balakrishnan algorithm requires the solution of an eigenvalue problem of size twice the order of the original system, but these are fairly small eigenvalue problems which can be solved efficiently and robustly using well-established factorization approaches. A structure-preserving algorithm for this task has been implemented as a FORTRAN subroutine for which we have used a MEX file to call it from MATLAB.  

Algorithm~\ref{alg1} is terminated in practice when the relative distance between $\omega_{r}$ and $\omega_{r-1}$ is less than a prescribed tolerance for some $r > 0$, or the number of iterations exceeds a specified integer. 
Formally, we terminate when
$$
	\abs{ \omega_r - \omega_{r-1} } < \varepsilon \cdot \frac{1}{2} \abs{\omega_r + \omega_{r-1}}
			 \quad \text{or} \quad r > r_{\max}.
$$
For our numerical experiments, we set $\varepsilon = 10^{-6}$ and $r_{\max} = 30$.

Algorithm~\ref{alg1} converges locally. To reduce the possibility of stagnating at a local maximizer
that is not a global maximizer, we initialize the algorithm with $r_0$ interpolation points $\omega_1,\,\dots,\,\omega_{r_0}$,
instead of only one. In our numerical experiments we have set $r_0 = 10$ as a default value, but there are more complicated examples that need a larger amount of initial interpolation points. For instance, the \texttt{peec} example {(see below for further details)} needs 80 initial points. We distribute the initial interpolation points equidistantly on the imaginary axis with the imaginary parts located in the interval $[0, \omega_{\max}]$, where $\omega_{\max}$ is a problem-dependent parameter that is highly influenced by the location of the poles of $H$. 

Another problem arises when the number of inputs and outputs is large. In this case, also the dimension $\tilde{n}$ of the middle factor $\widetilde{D}(s)$ of $\widetilde{H}(s)$ will grow with by $\min\{m,p\}$ in each interpolation step. To avoid a too fast growth of $\tilde{n}$ we have implemented an option in our implementation that allows to update the projection spaces only by using the singular vectors corresponding the the largest singular value of $H(\ri\omega_r)$. This means that in Algorithm \ref{alg1}, lines 11 and 13--17 disappear and line 12 is replaced by
\begin{align*}
 &\text{Compute the left and right singular vectors } w \text{ and } v \text{ of } \tilde{H}_{r-1}(\ri\omega_r) \text{ corresponding to } \\ 
 &\text{the largest singular value.} \\ 
 &\widetilde{V}_r \gets  D(\ri\omega_r )^{-1}B(\ri\omega_r)v\quad\text{and}\quad \widetilde{W}_r \gets \left( w^*C(\ri\omega_r)D(\ri \omega_r )^{-1} \right)^\ast.
\end{align*}
Similar changes are also made in lines 2--8. 
Note that in this way we may lose the Hermite interpolation property of the maximum singular values, since in general we only have $\sigma(\omega_k) \le \sigma_r(\omega_k)$, $k=1,\,\ldots,\,r$. Therefore, we also do not necessarily have local superlinear convergence. We have tested this option on the \texttt{mimo8x8\_system}, \texttt{mimo28x28\_system} and \texttt{mimo46x46\_system} examples, which have 8, 28, and 46 inputs and outputs, respectively.  The approach works well on these examples. However, a more rigorous analysis of this remains an open problem.  

In the next two subsections we report on the outcome of our numerical experiments. These have been performed on a machine with an 4 Intel\textsuperscript{\textregistered} Core\textsuperscript{\texttrademark} 3.30GHz i5-4590 CPUs and 16GB RAM in MATLAB 9.0.0.341360 (R2016a)
running on Linux version 3.12.67-64-default. First we test our algorithm on 33 linear systems taken from \cite{RomM06b,MarPR07,FreRM08,ChaVD02}
in Section \ref{sec:nexam_desc}. The data of these examples is freely available on the websites of Joost Rommes\footnote[2]{{see} \url{http://sites.google.com/site/rommes/software}} and the SLICOT benchmark collection\footnote[3]{{see} \url{http://slicot.org/20-site/126-benchmark-examples-for-model-reduction}}. The first 13 of these examples are standard state-space models ($E = I_n$), the other ones are descriptor systems with singular $E$. All these examples have transfer functions in $\cH_\infty^{p \times m}$, so in fact we compute the $\cH_\infty$-norm. Furthermore, we consider an example of a time-delay system provided in \cite{Gugercin2009}
in Section \ref{sec:nexam_delay}.

\subsection{Results for Descriptor Systems}\label{sec:nexam_desc}
In this subsection, we compare the results with the ones generated by the approach in \cite{Benner2014}, which is based on structured pseudospectra and locating their rightmost points in the complex plane repeatedly.
In this approach, perturbed transfer functions of the form 
$$
	H_{\Delta}(s) = C (sE - (A + B\Delta C))^{-1}B
$$
with $\Delta \in \C^{m \times p}$ are considered. There, a
perturbation $\Delta$ of minimal spectral norm such that 
the perturbed transfer function $H_{\Delta}$ is not in $\cH_\infty^{p \times m}$ is determined by a sequence of structured rank-1 perturbations.


\begin{sidewaystable}
 \caption{Numerical results for 33 test examples and comparison with the pseudospectral approach from \cite{Benner2014}} 
 \label{tab:results}
 \begin{tabular}{c|c|ccc|c|cc|cc|ccc}
  & & & & & & \multicolumn{2}{c}{computed $\cL_\infty$-norm} & \multicolumn{2}{c}{optimal frequency $\omega_\text{opt}$} & \multicolumn{3}{c}{time in s} \\
   \# & example & n & m & p & $n_{\rm iter}$ & \cite{Benner2014} & Algor.~\ref{alg1} & \cite{Benner2014} & Algor.~\ref{alg1} & \cite{Benner2014} & Algor.~\ref{alg1} & ratio \\
  \hline 
   1&  \texttt{build}  		& 48   &  1  & 1  &  6  &  5.27633e$-$03  & 5.27633e$-$03      & 5.20608e$+$00 & 5.20608e$+$00 &  1.06      & 0.08  & 14.0 \\
   2&  \texttt{pde}   		& 84   &  1  & 1  &  1  &  1.08358e$+$01  & 1.08358e$+$01      & 0.00000e$+$00 & 0.00000e$+$00 &  0.84      & 0.03  & 27.4 \\
   3&  \texttt{CDplayer}  	& 120   &  2  & 2 &  1	&  2.31982e$+$06  & 2.31982e$+$06      & 2.25682e$+$01 & 2.25682e$+$01 &  0.90      & 0.02  & 41.6 \\
   4&  \texttt{iss}   		& 270  &  3  & 3  &  7	&  1.15887e$-$01  & 1.15887e$-$01      & 7.75093e$-$01 & 7.75093e$-$01 &  0.85      & 0.24  & 3.5 \\
   5&   \texttt{beam} 		& 348   &  1  & 1  & 1  &  4.55487e$+$03  & 4.55487e$+$03      & 1.04575e$-$01 & 1.04575e$-$01 & 11.06      & 0.08 & 135.8 \\
   6&  \texttt{S10PI\_n1} 	& 528   & 1   & 1  & 7	&  3.97454e$+$00  & 3.97454e$+$00      & 7.53151e$+$03 & 7.53151e$+$03 &  0.79      &  0.08  & 10.3 \\
   7&  \texttt{S20PI\_n1}   	&  1028  &   1 &  1 & 5	&  3.44317e$+$00  & 3.44317e$+$00      & 7.61831e$+$03 & 7.61831e$+$03 &  1.79      & 0.07  & 24.0 \\
   8&  \texttt{S40PI\_n1}    	& 2028   &  1  & 1  & 7	&  3.34732e$+$00  & 3.34732e$+$00      & 6.95875e$+$03 & 6.95875e$+$03 &  1.95      & 0.14  & 13.6 \\
   9&  \texttt{S80PI\_n1}   	&  4028  &  1  &  1 & 5	&  3.37016e$+$00  & 3.37016e$+$00      & 6.96149e$+$03 & 6.96149e$+$03 &  3.84      & 0.19  & 20.3 \\
   10&  \texttt{M10PI\_n1}   	&  528  &  3  &  3  & 7 &  4.05662e$+$00  & 4.05662e$+$00      & 7.53181e$+$03 & 7.53181e$+$03 &  1.21     & 0.35  & 3.5 \\
   11&  \texttt{M20PI\_n1}   	&  1028  &  3  &  3 & 12&  3.87260e$+$00  & 3.87260e$+$00      & 5.06412e$+$03 & 5.06412e$+$03 &  1.42     & 0.85  & 1.7 \\
   12&  \texttt{M40PI\_n1}    	&  2028  &  3  & 3 & 8	&  3.81767e$+$00  & 3.81767e$+$00      & 5.07107e$+$03 & 5.07107e$+$03 &  2.24      & 0.51  & 4.4 \\
   13&   \texttt{M80PI\_n1}   	&  4028  &  3  &  3 & 9	&  3.80375e$+$00  & 3.80375e$+$00      & 5.07279e$+$03 & 5.07279e$+$03 &  3.82     & 0.83  & 4.6 \\ \hline
   14&   \texttt{peec}   	&  480  &  1  & 1  & 1	&  3.52624e$-$01  & 3.52610e$-$01      & 5.46349e$+$00 & 5.46349e$+$00 &  9.26      & 2.13 & 4.3 \\
   15&   \texttt{S10PI\_n}     	&  682  &  1  &  1 & 7  &  3.97454e$+$00  & 3.97454e$+$00      & 7.53151e$+$03 & 7.53151e$+$03 &  1.03     & 0.09  & 11.3 \\
   16&   \texttt{S20PI\_n}  	&  1182  & 1  & 1  & 5  &  3.44317e$+$00  & 3.44317e$+$00      & 7.61831e$+$03 & 7.61831e$+$03 &  1.90     & 0.08  & 24.7 \\
   17&   \texttt{S40PI\_n}   	&  2182  & 1  &  1 & 7	&  3.34732e$+$00  & 3.34732e$+$00      & 6.95875e$+$03 & 6.95875e$+$03 &  2.12      & 0.15 & 14.6 \\
   18&   \texttt{S80PI\_n}    	& 4182   & 1  & 1  & 5	&  3.37016e$+$00  & 3.37016e$+$00      & 6.96149e$+$03 & 6.96149e$+$03 &  3.96     & 0.20  & 20.1 \\
   19&   \texttt{M10PI\_n}   	&  682  &  3  &  3 & 7	&  4.05662e$+$00  & 4.05662e$+$00      & 7.53181e$+$03 & 7.53181e$+$03 &  1.40     & 0.35  & 4.0 \\
   20&   \texttt{M20PI\_n}    	& 1182   & 3  & 3  & 10	&  3.87260e$+$00  & 3.87260e$+$00      & 5.06412e$+$03 & 5.06412e$+$03 &  1.44     & 0.60  & 2.4 \\
   21&   \texttt{M40PI\_n}    	&  2182  &  3  & 3 &  8	&  3.81767e$+$00  & 3.81767e$+$00      & 5.07107e$+$03 & 5.07107e$+$03 &  2.12     & 0.51  & 4.1 \\
   22&   \texttt{M80PI\_n}   	&  4182  &  3  & 3 & 9 	&  3.80375e$+$00  & 3.80375e$+$00      & 5.07279e$+$03 & 5.07279e$+$03 &  3.96      & 0.85  & 4.7 \\
   23&   \texttt{bips98\_606}   & 7135   &  4  &  4 & 1	&  2.01956e$+$02  & 2.01956e$+$02      & 3.81763e$+$00 & 3.81762e$+$00 & 14.18    & 0.66  & 21.4 \\
   24&   \texttt{bips98\_1142}  &  9735  &  4  &  4  & 1 &  1.60427e$+$02  & 1.60427e$+$02      & 4.93005e$+$00 & 4.93006e$+$00 & 29.25    & 0.83  & 35.4 \\
   25&   \texttt{bips98\_1450}   &  11305  &  4  & 4 & 1 &  1.97389e$+$02  & 1.97389e$+$02      & 5.64575e$+$00 & 5.64571e$+$00 & 26.16     & 0.94  & 27.9 \\
   26&   \texttt{bips07\_1693}   &  13275  &  4  & 4 & 1 & 2.04168e$+$02  & 2.04168e$+$02      & 5.53766e$+$00 & 5.53765e$+$00 & 66.59      & 1.07 & 62.3 \\
   27&   \texttt{bips07\_1998} 	& 15066   &  4  & 4  & 2 &  1.97064e$+$02  & 1.97064e$+$02      & 6.39968e$+$00 & 6.39960e$+$00 & 40.37     & 1.63  & 24.8 \\
   28&   \texttt{bips07\_2476}   & 16861   &  4  &  4 &	2 &  1.89579e$+$02  & 1.89579e$+$02      & 5.88971e$+$00 & 5.88973e$+$00 & 64.88     & 1.96  & 33.1 \\
   29&   \texttt{bips07\_3078}   &  21128  &  4  &  4 & 1 & 2.09445e$+$02  & 2.09445e$+$02      & 5.55792e$+$00 & 5.55793e$+$00 &  35.18      & 2.08  & 16.9 \\
   30&   \texttt{xingo\_afonso\_itaipu} &  13250  &  1  & 1 & 2	&  4.05605e$+$00    & 4.05605e$+$00   & 1.09165e$+$00  & 1.09165e$+$00  &  14.38  & 0.56 & 24.6 \\
   31&   \texttt{mimo8x8\_system}    &  13309  &  8  &  8  & 2 &  5.34292e$-$02    & 5.34292e$-$02   & 1.03313e$+$00  & 1.03312e$+$00  &  26.74 & 1.27  & 21.0 \\
   32&   \texttt{mimo28x28\_system}  &  13251  & 28   & 28  & 3	&  1.18618e$-$01    & 1.18618e$-$01   & 1.07935e$+$00  & 1.07935e$+$00  &  24.78 & 2.62  & 9.5 \\
   33&   \texttt{mimo46x46\_system}  & 13250   & 46   &  46  & 3 &  2.05631e$+$02    & 2.05631e$+$02   & 1.07908e$+$00  & 1.07908e$+$00  & 36.84 & 3.76  & 9.8 \\
   \end{tabular}
 \end{sidewaystable}

Table \ref{tab:results} summarizes the results of the 33 numerical experiments. For all examples, the correct norm value has been found up the termination tolerance. In this table, the number of additional iterations after the construction of the initial reduced function needed to retrieve the $\cL_\infty$-norm by Algorithm~\ref{alg1} up to the prescribed relative tolerance $\varepsilon = 10^{-6}$ is denoted by $n_{\rm iter}$. The order of the system, the input dimension, and the output dimension are denoted by $n,\,m,\,p$, respectively. It is evident from Table~\ref{tab:results} that the correct value of the $\cH_\infty$-norm is found by Algorithm~\ref{alg1} for each of the problems. In terms of the runtime, Algorithm~\ref{alg1} outperforms the pseudospectral approach. The ratios between the time required by the pseudospectral approach and that 
required by Algorithm~\ref{alg1} are listed in the last column of Table~\ref{tab:results}.

Finally, local superlinear convergence consistent with Theorem~\ref{thm:quad_conv} is observed in all cases. 
Specifically, for the \texttt{S80PI\_n} example, the errors of the iterates are reported in Table~\ref{tab:quad}.
Five additional iterations after the construction of the initial reduced function suffice to compute the $\cL_\infty$-norm with a desired relative tolerance of $\varepsilon = 10^{-6}$. 
In fact, we see that once the algorithm started converging to a local maximizer, it only needs one or two more iterations until convergence.    
\begin{table}[t]
\caption{The errors of the iterates of Algorithm~\ref{alg1}, the ratios of the errors of the iterates, and the errors of the largest singular values
at these iterates are listed for the \texttt{S80PI\_n} example. Here, the short-hands 
$\sigma_r := \sigma_r(\omega_{r+1})$ and $\sigma_\ast := \sigma(\omega_\ast)$ are used. As the ``exact'' solution we have taken the one we obtain after iteration 6.}
\label{tab:quad}
\begin{center}
\begin{tabular}{l|ccc}
	Iteration $\#$ ($r$)	&	$|\omega_{r+1} - \omega_\ast|$	 &	$\; |\omega_{r} - \omega_\ast| / |\omega_{r-1} - \omega_\ast| \;$ & $|\sigma_{r} - \sigma_\ast | $	\\	
\hline	
   1 (Initial model)  &  6.718e$+$02      &  ---			& 1.207e$+$02  \\
   2                  &  2.595e$+$03      &  3.863e$+$00	& 4.388e$+$00   \\
   3                  &  1.156e$+$01     &  4.455e$-$03	& 2.905e$+$00    \\
   4                  &  6.571e$-$01       &  5.684e$-$02	        & 8.662e$-$03    \\
   5                  &     0 	    &          0   	& 8.878e$-$09	\\
\end{tabular}
\end{center}
\end{table}

\begin{figure}[ht!]
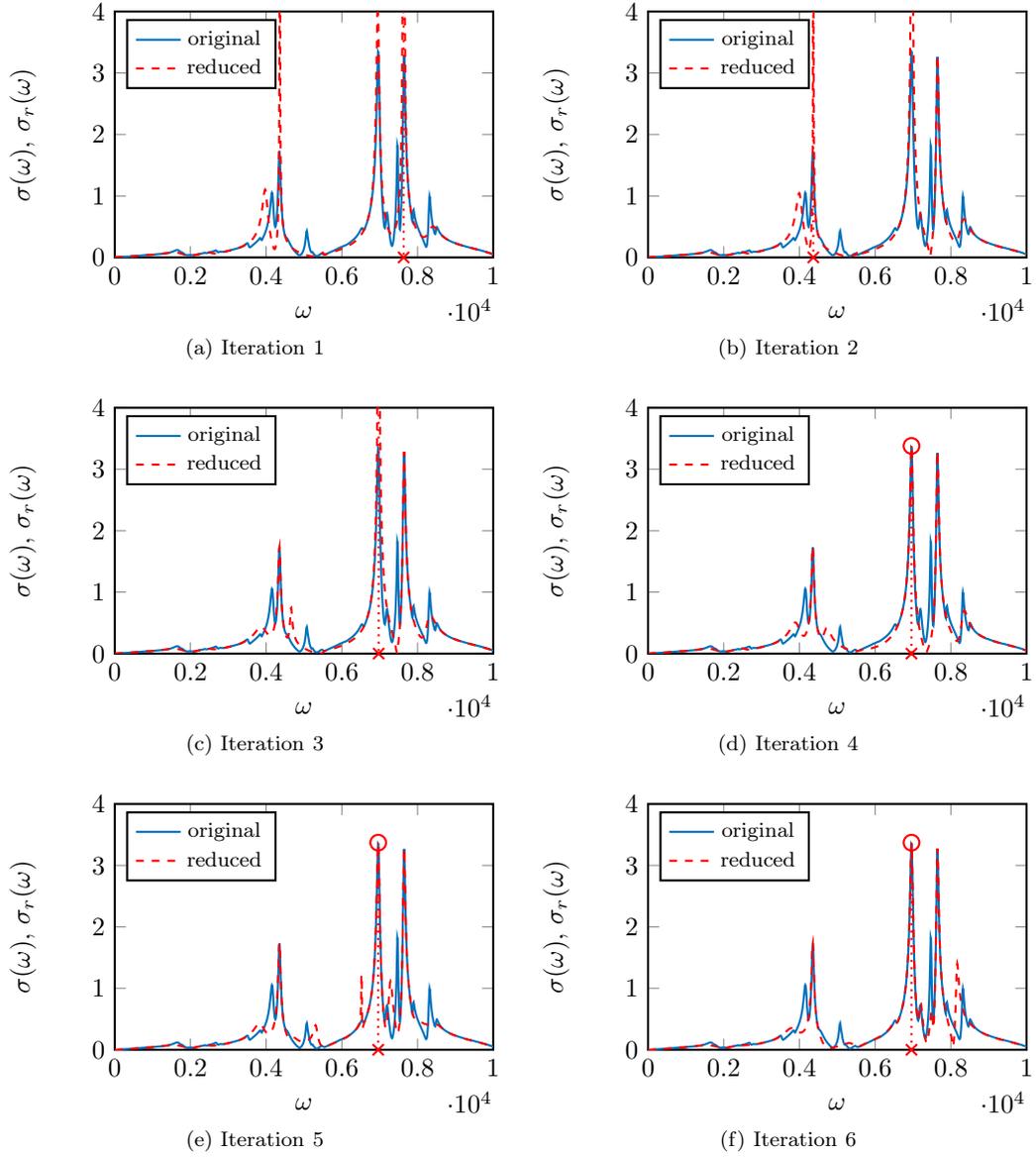

 \centering
 \subfloat[][Iteration 1]{\input{initial}} \quad
 \subfloat[][Iteration 2]{\input{iter1}} \\
 \subfloat[][Iteration 3]{\input{iter2}} \quad
 \subfloat[][Iteration 4]{\input{iter3}} \\
 \subfloat[][Iteration 5]{\input{iter4}} \quad
 \subfloat[][Iteration 6]{\input{iter5}}
 \caption{Intermediate reduced functions obtained by Algorithm \ref{alg1} for the \texttt{S80PI\_n} example. The original function is depicted in blue, while the reduced functions are represented by the dashed red lines. The red crosses and circles indicate the locations of the maximizers and the $\cL_\infty$-norms of the reduced functions, respectively. In the first iteration (on {the} top left), the reduced model represents the initial model obtained from 10 equally spaced interpolation points.} 
 \label{fig:iters}
\end{figure}

However, in some numerical examples we also observe that many iterations may be needed until convergence to a local maximizer takes place. If there exist many local maximizers of $\sigma( \cdot)$, then the algorithms often collects more global information of $H$ in the beginning and only starts converging to a local maximizer after a certain number of iterations. In particular, this is the case in examples \# 6--13 and \# 15--22. An illustration of this fact is given in Figure \ref{fig:iters}, where the intermediate reduced functions for the \texttt{S80PI\_n} example are depicted.

\begin{figure}
\begin{center}
\begin{tikzpicture}
\definecolor{mycolor1}{rgb}{0.00000,0.44700,0.74100}%
\begin{axis}[%
width=4.0in,
height=2.0in,
line width = 0.8pt,
at={(0.75in,0.478in)},
scale only axis,
xmin=1,
xmax=30,
xlabel={\texttt{ninit}},
ymin=0,
ymax=25,
ylabel={\large $n_{\rm iter}$},
legend pos=north east,
axis y line = left,
]
\addplot [color=mycolor1,solid]
  coordinates{%
  (1,19)
  (2,20)
  (3,16)
  (4,15)
  (5,6)
  (6,14)
  (7,10)
  (8,11)
  (9,5)
  (10,9)
  (11,11)
  (12,10)
  (13,3)
  (14,6)
  (15,5)
  (16,5)
  (17,3)
  (18,2)
  (19,1)
  (20,1)
  (21,1)
  (22,1)
  (23,1)
  (24,1)
  (25,1)
  (26,3)
  (27,1)
  (28,1)
  (29,1)
  (30,1)
};\label{plot_one}
\addlegendentry{\footnotesize{$n_{\rm iter}$}}
\end{axis}
\begin{axis}[%
width=4.0in,
height=2.0in,
line width = 0.8pt,
at={(0.75in,0.478in)},
scale only axis,
ymin=0,
ymax=1.5,
xmin=1,
xmax=30,
ylabel={time in seconds},
axis x line=none,
legend pos=north east,
axis y line = right,
]
\addlegendimage{/pgfplots/refstyle=plot_one}\addlegendentry{\footnotesize{$n_{\rm iter}$}}
\addplot [color=red,dashed]
  coordinates{%
  (1,1.19682)
  (2,1.43705)
  (3,0.995531)
  (4,1.01999)
  (5,0.293387)
  (6,1.06433)
  (7,0.681273)
  (8,0.891439)
  (9,0.40414)
  (10,0.840512)
  (11,1.28615)
  (12,1.19388)
  (13,0.421878)
  (14,0.807917)
  (15,0.786009)
  (16,0.834236)
  (17,0.654967)
  (18,0.534017)
  (19,0.450442)
  (20,0.480074)
  (21,0.540249)
  (22,0.536108)
  (23,0.591584)
  (24,0.613634)
  (25,0.65463)
  (26,1.16573)
  (27,0.717959)
  (28,0.805636)
  (29,0.835097)
  (30,0.796568)
};
\addlegendentry{\footnotesize{runtime}}
\end{axis}
\end{tikzpicture}
\caption{Behavior of Algorithm~\ref{alg1} {on the \texttt{M80PI\_n}} example with respect to the number of initial interpolation points. The value $n_{\rm iter}$ refers to the number of additional iterations after the construction of the initial reduced function.}
\label{fig:interpolpoints}
\end{center}
\end{figure}
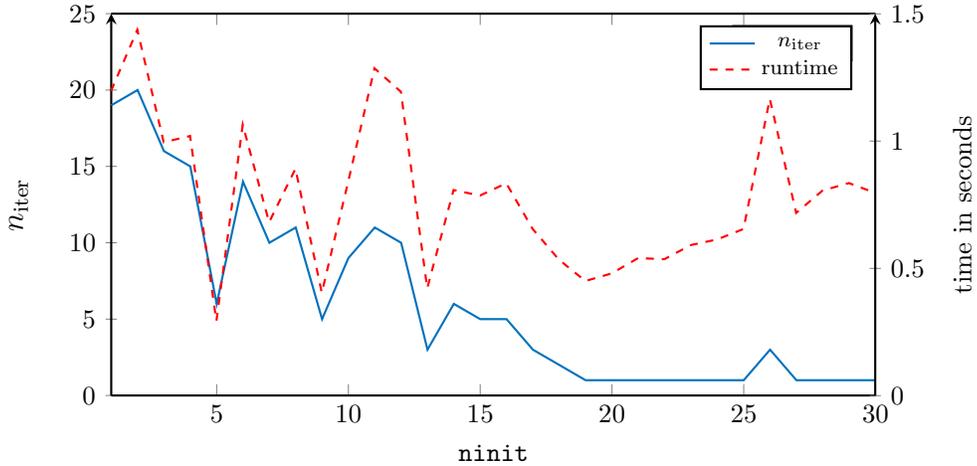

A further strong influence of the performance of Algorithm~\ref{alg1} is the number and location of the initial interpolation points. In the examples above, we have usually taken 10 initial points distributed equidistantly in an interval $[0,\omega_{\max}]$. Often, the algorithm also converges to the correct global maximizer with fewer initial points but then it may happen that more iterations are needed, since less global information is known. To illustrate this behavior we consider the \texttt{M80PI\_n} example. For this example we generate interpolation points by the MATLAB command
\begin{verbatim}
options.initialPoints = linspace( 0.1, 10000, ninit );
\end{verbatim}
and we let \texttt{ninit} grow from 1 to 30. For all 30 configurations, the correct norm value has been computed. The results are depicted in Figure~\ref{fig:interpolpoints}. It can be seen that there is a certain trade-off between the number of initial interpolation points; a larger value of \texttt{ninit} may drastically reduce the number of additional iterations, but it also increases the subspace dimensions, which results in more effort for solving the small intermediate problems. For our example, values of about 20 initial interpolation points result in the best behavior (except for some smaller values, where the global optimizer has already been almost hit).

\subsection{Results for Time-Delay Systems}\label{sec:nexam_delay}

Next we test our approach on transfer functions of time-delay systems. Our experiments are performed on the following example taken from \cite{Gugercin2009}.


\begin{example}
\label{ex:delay}
Consider the delay system
\begin{equation}\label{eq:delay_sys_onedelay}
	Ex'(t) = A_0 x(t) +  A_1 x(t - \tau) + B u(t), \quad
	  y(t) = C x(t)
\end{equation}
with $E = \theta I_n + T$, $A_0 = \frac{1}{\tau} \left( \frac{1}{\beta}+1\right)(T - \theta I_n)$, $A_1 = \frac{1}{\tau} \left( \frac{1}{\beta} - 1\right)(T - \theta I_n)$, where $T$ is the $n \times n$ matrix with ones on the subdiagonal and superdiagonal, as well as in the entries at position $(1,1)$ and $(n,n)$, and zeros elsewhere. The scalars $\beta$ and $\theta$ are parameters, and $\tau$ is the delay parameter.

We choose $\tau = 1$, $\beta = 0.01$, and $\theta = 5$.
Additionally, we set $B = e_1 + e_2$, the sum of the first two columns of the $n\times n$ 
identity matrix, $C = B^*$, and experiment with various values of $n$.
\end{example}

Since Example \ref{ex:delay} has the non-rational transfer function
\begin{equation*}
 H(s) = C\left(sE - A_0 - \mathrm{e}^{-s\tau}A_1 \right)^{-1} B,
\end{equation*}
the Boyd-Balakrishnan algorithm cannot be applied here. Instead we use \texttt{eigopt} \cite{Mengi2014} to solve the small subproblems in Algorithm~\ref{alg1}. The Matlab package \texttt{eigopt} requires additional inputs. Specifically, a frequency interval in which the $\cL_\infty$-norm is attained has to be supplied. We choose the interval $[0, 50]$, in which $\omega_*$ is located. In this interval, 8 local maxima of $\sigma(\cdot)$ {can be found}. Outside this interval, there exist infinitely many more such local maxima, but they result in much smaller maximum singular values. A second parameter the user has to supply is a global lower bound $\gamma$ on the second derivative of $-\sigma(\cdot)$. In our example, the minimum of this second derivative is always about $-93.08$, so we choose $\gamma = -100$. The value of $\gamma$ has a strong influence on the runtime; the lower $\gamma$, the more piecewise quadratic support functions are constructed by \texttt{eigopt} which increases the computational complexity. 

The runtimes and the runtime ratios between \texttt{eigopt} and Algorithm~\ref{alg1} are given in Table~\ref{tab:delay}. For all values of $n$, \texttt{eigopt} and Algorithm~\ref{alg1} return the same (correct) value of the $\cL_\infty$-norm, namely $\left\| H \right\|_{\cL_\infty} = 0.23766$. This value is attained for $\omega_* = 3.07547$. After the construction of the initial reduced transfer function, Algorithm~\ref{alg1} only needs one more iteration until convergence. The table also shows that Algorithm~\ref{alg1} is only more efficient for larger values of $n$. This is because the computation of $H$ and its singular values becomes a dominant factor for larger $n$. For smaller values of $n$, Algorithm~\ref{alg1} carries out two calls of \texttt{eigopt}; for this reason, it needs almost the double the time to solve the original problem by a single run of \texttt{eigopt}. 
 
%
%
%
%
%
%

\begin{table}[tb]
\begin{center}
\caption{Comparison of \texttt{eigopt} and Algorithm \ref{alg1} on Example \ref{ex:delay}}
\label{tab:delay}
\begin{tabular}{c|cc|c} 
    &    \multicolumn{2}{c}{time in seconds} &   \\  
  $n$ &   \texttt{eigopt} & Algor.~\ref{alg1} & ratio  \\  
\hline 
100  & 0.92 & 1.83 & 0.51 \\
300  & 0.98 & 1.86 & 0.53 \\
1000  & 1.25 & 1.85 & 0.67 \\
3000  &  1.76 & 1.86 & 0.95 \\
10000  & 3.92 & 1.89 & 2.07 \\
30000  & 10.47 & 2.02 & 5.19 \\
100000 & 36.36 & 2.52 & 14.38 \\
300000 & 113.10  & 3.80 & 29.75 \\
1000000 & 403.13 & 9.09 & 44.34
\end{tabular}
\end{center}
\end{table}

\subsection{Limitations of the Method}
\begin{figure}[tb]
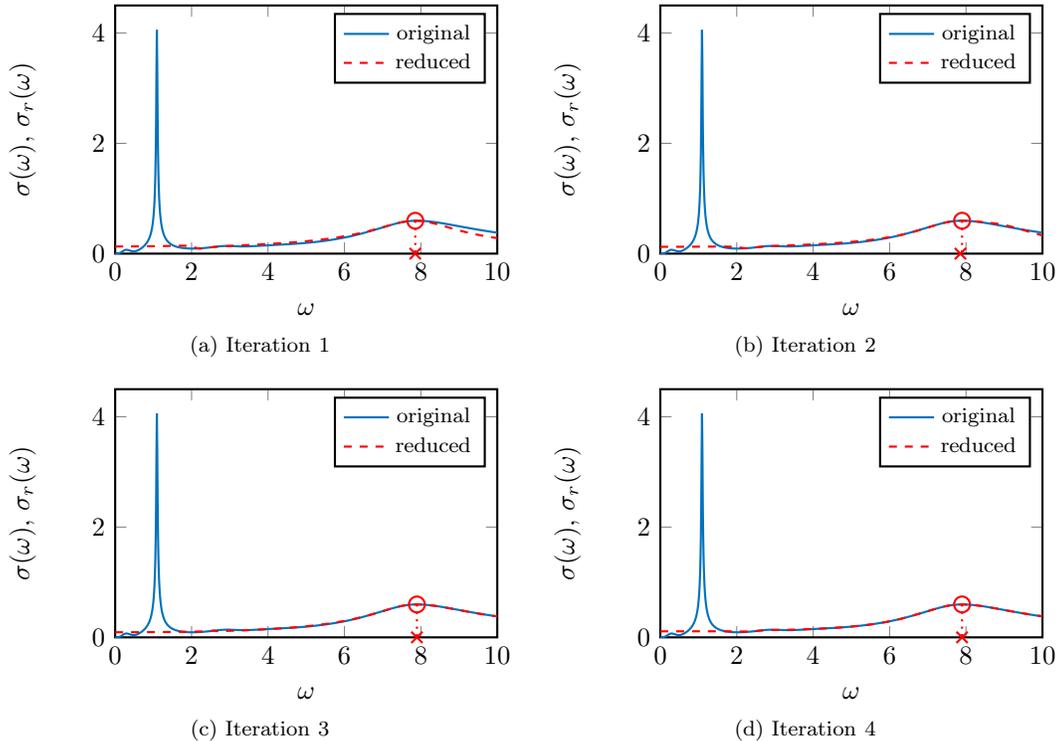

 \centering
 \subfloat[][Iteration 1]{\input{initial_xingo}} \quad
 \subfloat[][Iteration 2]{\input{iter1_xingo}} \\
 \subfloat[][Iteration 3]{\input{iter2_xingo}} \quad
 \subfloat[][Iteration 4]{\input{iter3_xingo}} 
 \caption{Intermediate reduced functions obtained by Algorithm \ref{alg1} for the \texttt{xingo\_afonso\_itaipu} example. The original function is depicted in blue, while the reduced functions are represented by the dashed red lines. The red crosses and circles indicate the locations of the maximizers and the $\cL_\infty$-norms of the reduced functions, respectively. The initial reduced model (on {the} top left) is obtained from 2 interpolation points at $\omega = 2.5$
 and $\omega = 7.5$.}
 \label{fig:-(}
\end{figure}

As mentioned above, our algorithm converges only locally. (The same property holds for all other methods for large-scale $\cL_\infty$-norm calculationsto this date.) It is important to interpolate $H$ at the parts of the imaginary axis that are close to the poles of $H$. If not enough interpolation points are taken, then the global maximizer of $\sigma(\cdot)$ may be missed. To illustrate this, consider the \texttt{xingo\_afonso\_itaipu} example but with only two initial interpolation points $2.5\ri$ and $7.5\ri$. With these points only, the global maximizer at 1.092 is not detected, instead the algorithm converges to the local maximizer at 7.897 which is not a global maximizer. The intermediate iterates are depicted in Figure \ref{fig:-(}. 
A remedy to this problem is to use more initial interpolation points for the initial iteration. 

\section{Concluding Remarks}\label{sec:conclude}
We have introduced an approach for the computation of the $\cL_\infty$-norm of an $\cL_\infty$-function 
of the form $H(s) = C(s)D(s)^{-1}B(s)$ in the large-scale setting, 
{i.\,e.}, the middle factor is the inverse of a large-scale meromorphic matrix-valued function, and $C(s),\, B(s)$ are meromorphic functions mapping to short-and-fat and tall-and-skinny matrices, respectively.
Our approach is based on a subspace projection idea that is frequently used in model order reduction. More precisely, we approximate the given $\cL_\infty$-function by a reduced function obtained by employing two-sided projections on the factors of the original $\cL_\infty$-function. The middle factor of the resulting reduced function is of much smaller dimension. We compute the $\cL_\infty$-norm of the reduced function by established methods. Then we expand the projection spaces by using the singular vectors of the original function at the point on imaginary axis, where the $\cL_\infty$-norm of the reduced function is attained.
We have proven that our selection strategy for the subspaces leads to Hermite interpolation properties between the largest singular values of the original and reduced functions. These Hermite interpolation properties in turn give rise to a superlinear convergence with respect to the subspace dimension.

We have demonstrated on various numerical examples that our method can lead to substantial speedups compared to known methods. Moreover, it can be applied to a much larger class of functions such as transfer functions of delay systems. Thus, our method may lead to significant computational benefits in the field of $\cH_\infty$-optimization.

\vskip 2ex

\noindent
\textbf{Acknowledgements.} The authors are grateful to two anonymous reviewers and Daniel Szyld for their invaluable feedback.



\bibliography{large_hinfinity_new}

\end{document}